\documentclass[12pt]{article}
\usepackage{a4wide}
\usepackage{amssymb}
\usepackage{amsfonts}
\usepackage{amsmath}

\date{}

\newtheorem{proposition}{Proposition}[section]
\newtheorem{theorem}[proposition]{Theorem}
\newtheorem{lemma}[proposition]{Lemma}

\newtheorem{corollary}[proposition]{Corollary}

\def\GK{{\rm  GK}\,}

\def\der{\partial }

\def\nFM0{{\nu }_{F,M_0}}
\def\nFN0{{\nu }_{F,N_0}}
\def\nGN0{{\nu }_{G,N_0}}

\def\N0{ {\bf N}_0 }

\def\t{\otimes}
\def\g{\gamma}

\def\ra{\rightarrow}
\def\lra{\leftrightarrow}
\def\Xpm{X^{\pm }}

\def\s{\sigma}

\def\l1{{\lambda}_1}

\def\m{{\bf m}}
\def\a{\alpha}
\def\a0{ {\alpha }_0}
\def\a1{ {\alpha }_1}

\def\l{\lambda}
\def\o{\omega}

\def\nFGM0{{\nu }_{F,G,M_0}}

\def\nFN0{{\nu}_{F,N_0}}

\def\sm{{\sigma}^m}
\def\si{{\sigma}^i}

\def\sm1{{\sigma}^{-1}}

\def\smtp1{{\sigma}^{-t+1}}

\def\o{\omega }
\def\S1{S^{-1}}

\def\Xpm1{X^{\pm 1}_1}

\def\sPM1{{\sigma }^{\pm 1}}
\def\sMP1{{\sigma }^{\mp 1 }}

\def\d{\delta}

\def\L{\Lambda}

\def\G{\Gamma}

\def\OO{{\cal O}}
\def\CA{{\cal A}}

\def\CD{{\cal D}}

\def\Ytm1{Y^{t-1}}
\def\Yim1{Y^{i-1}}

\def\CL{{\cal L}}

\def\CG{{\cal G}}
\def\CH{{\cal H}}

\def\supp{{\rm supp }}

\def\Aut{{\rm Aut}}
\def\bK{\overline{K}}

\def\bA{\overline{A}}
\def\Der{{\rm Der }}
\def\ad{{\rm ad }}
\def\dim{{\rm dim }}

\def\ker{ {\rm ker } }

\def\gr{ {\rm gr} }

\def\D{ \Delta }

\def\Ev{ {\rm Ev} }
\def\CR{ {\cal R} }

\def\h'{ \tilde{h} }

\def\End{ {\rm End} }
\def\CB{ {\cal B} }
\def\CC{ {\cal C} }
\def\CE{ {\cal E} }
\def\ord{{\rm ord}}

\begin{document}

\author{V. V.\  Bavula }

\title{Dixmier's Problem 6 for the Weyl Algebra (the Generic Type Problem)}

\maketitle
\begin{abstract}
In the fundamental paper \cite{Dix}, J. Dixmier posed six problems
for the Weyl algebra $A_1$ over a field $K$ of characteristic
zero. Problem 3 was solved by Joseph and Stein \cite{josclA1}
(using results of McConnel and Robson \cite{MRext}), problem 5 was
solved by the author in \cite{BavDP5}. Using a (difficult)
polarization theorem for the Weyl algebra $A_1$  Joseph
\cite{josclA1} solved problem 6. Problems 1, 2, and 4 are still
open. Note that these problems make sense for non-commutative
algebras of Gelfand-Kirillov $<3$, and some of them (after minor
modifications) make sense for an arbitrary noncommutative algebra.

In this paper a short proof is given to Dixmier's problem $6$ for
 many noncommutative  algebras $A$ of Gelfand-Kirillov $<3$ (a typical example is
 the ring of differential operators $\CD (X)$ on a smooth irreducible
  algebraic curve $X$).  An
  affirmative answer to this problem  leads to clarification of
  the structure of maximal commutative subalgebras of these
  algebras $A$ (a typical example of the algebra $A$ is
 any  noncommutative subalgebra   of Gelfand-Kirillov $<3$
  of the division ring $Q(\CD (X))$ for the algebra $\CD (X)$),
    and the result is rather surprising:
  for a given maximal commutative subalgebra $C$ of the algebra
  $A$,  (almost) all non-central elements of it have the {\em same}
  type, more precisely, have exactly one of the following types: (i) strongly nilpotent, (ii)
  weakly nilpotent, (iii) generic, (iv) generic except for a
  subset $K^*a+Z(A)$  of strongly semi-simple elements,
  (iv) generic except for a subset $K^*a+Z(A)$  of weakly semi-simple
  elements, where $K^*:=K\backslash \{ 0\} $
  and $Z(A)$ is the centre of the algebra $A$.

  For an arbitrary algebra $A$, Dixmier's problem  6 is essentially
  a question: {\em whether an inner derivation of the algebra $A$
  of the type $\ad \, f(a)$, $a\in A$, $f(t)\in K[t]$, $\deg_t
  (f(t))>1$, has a nonzero eigenvalue.} We prove that the answer
  is negative for many classes of algebras ({\em eg}, rings of
  differential operators $\CD (Y)$ on smooth irreducible algebraic
   varieties, all prime factor algebras of the universal enveloping
   algebra $U(\CG )$  of a completely solvable algebraic Lie algebra $\CG $).

\end{abstract}


\section{Introduction}

In this paper,   $K$ is a field of characteristic zero and $\bK $
is its algebraic closure.

The {\em (first) Weyl algebra} $A_1$
 is an associative algebra generated over the field $K$ by elements $x$ and $\der $
 that satisfy the defining relation $\der x-x\der =1$. The Weyl
algebra $A_1$ is a central, simple, Noetherian domain of
Gelfand-Kirillov dimension $2$ which is canonically isomorphic to
the ring of differential operators $K[x][ { d\over dx}]$ $(x\lra
x$, $\der \lra { d\over dx})$ with coefficients from the
polynomial algebra $K[x]$. The $n$'th {\em Weyl algebra} $A_n$ is
the tensor product $A_1\t \cdots \t A_1$ of $n$ copies of the
first Weyl algebra.

In the seminal paper \cite{Dix}, Dixmier initiated a systematic
study of the structure of the first Weyl algebra. At the end of
this paper he posed six problems (questions).  Problem 3 has been
solved by Joseph and Stein \cite{josclA1} (using results of
McConnel and Robson \cite{MRext}), problem 5 has been solved by
the author in \cite{BavDP5}. Using a (difficult) polarization
theorem for the Weyl algebra $A_1$ Joseph \cite{josclA1} solved
problem 6. Problems 1, 2, and 4 are still open. {\bf  Dixmier's
problem 1}: {\em whether an algebra endomorphism of the first Weyl
algebra is an algebra automorphism?}  He writes in the paper on
page 242 that ``{\em A. A. Kirillov informed me that the Moscow
school also considered this problem.}''  A positive answer to a
similar
 problem but for the $n$'th Weyl algebra implies the {\em Jacobian
 Conjecture} (see the paper of  Bass, Connel and  Wright
 \cite{BCW} for details, and  \cite{Bavrenqu} on a progress on these two problems).

 Let $A$ be an algebra over the field
$K$, and let $\Der_K(A)$ be the set of all $K$-derivations of the
algebra $A$. We can identify the algebra $A$ with its image in the
algebra $\bA :=\bK \t_K A$ via the $K$-algebra monomorphism $A\ra
\bA $ , $a\ra  1\t a$. Then $\Der_K(A) \subseteq \Der_{\bK }(\bA
)$. $\Der_K(A)$ is a Lie algebra with the bracket given by the
commutator of derivations $([\d ,\der ]:=\d \der -\der \d )$. For
each derivation  $\d \in \Der_K(A)$ of the algebra $A$ one can
attach four subalgebras of $A$: the {\em kernel}
  $C(\d ):=\ker \, \d $; the ($\mathbb{N}$-filtered) {\em nil-algebra}
$N(\d ):=\bigcup_{i\geq 0} N(\d , i)$  where $N(\d , i):= \ker \,
\d^{i+1}$;  the {\em eigenvalue algebra} $D(\d ):=A\cap D(\d , \bA
)$ for $\d $ where  the subalgebra $D(\d , \bA)$ of $\bA $ is  the
direct sum
 $$ \bigoplus_{\l \in \Ev (\d , \bA )} D(\d , \l , \bA )$$ of the
eigenspaces  $D(\d , \l ,\bA ):=\{ a\in \bA \, |\, \d (a) =\l a\}$
belonging to the eigenvalue $\l \in \bK $,  $\Ev (\d , \bA )$ is
the set of all eigenvalues for $\d $ in $\bA $; the {\em torsion
algebra} $F(\d ):=\{ a\in A\, | \, \dim_K \, K[\d ]a<\infty \}$
for $\d $.  Clearly, $C(\d )=N(\d )\cap D(\d )$, $N(\d )\subseteq
F(\d )$, and $D(\d )\subseteq F(\d )$.

The set $\Der_K(A)\backslash  \{ 0\} $ is a disjoint union of its
eleven subsets
$$ \Der_K(A)\backslash  \{ 0\} =\bigcup_{i=1}^{11}\, \widetilde{\Delta}_i(A)  $$
where
\begin{enumerate}
\item $\widetilde{\Delta}_1(A)=\{ \d :\; N(\d )=A,\;
D(\d )=C(\d )\}$, a set of {\em strongly nilpotent} derivations.
\item $\widetilde{\Delta}_2(A)=\{ \d :\; C(\d ) \neq  N(\d )\neq A,\;
 D(\d )=C(\d )\}$, a set of {\em weakly nilpotent} derivations.
\item $\widetilde{\Delta}_3(A)=\{ \d :\; N(\d )=C(\d ), \;  D(\d )=A
\}$, a set of {\em strongly semi-simple} derivations.
\item $\widetilde{\Delta}_4(A)=\{ \d :\;  N(\d )=C(\d ),\;  C(\d )\neq D(\d )=F(\d )\neq A\}$, a set of
{\em weakly semi-simple} derivations.
\item $\widetilde{\Delta}_5(A)=\{ \d :\; N(\d )=D(\d )=C(\d )\}$, a set of {\em generic} derivations.

\item $\widetilde{\Delta}_6(A)=\{ \d :\; N(\d )=C(\d )\neq D(\d )\neq F(\d )=A\}$, a set of
{\em strongly jordan} derivations.
\item $\widetilde{\Delta}_7(A)=\{ \d :\; N(\d )=C(\d )\neq D(\d )\neq F(\d )\neq A\}$, a set of
{\em weakly jordan} derivations.
\item $\widetilde{\Delta}_8(A)=\{ \d :\; N(\d )\neq C(\d )\neq D(\d ),\;  N(\d )+D(\d )\neq F(\d )=A\}$, a set of
{\em strongly nil-jordan} derivations.
\item $\widetilde{\Delta}_9(A)=\{ \d :\; N(\d )\neq C(\d )\neq D(\d ),\;  N(\d )+D(\d )\neq F(\d )\neq A\}$, a set of
{\em weakly nil-jordan} derivations.
\item $\widetilde{\Delta}_{10}(A)=\{ \d :\; N(\d )\neq C(\d )\neq D(\d ), N(\d )+D(\d )=A\}$, a set of
{\em strongly nil-semi-simple} derivations.
\item $\widetilde{\Delta}_{11}(A)=\{ \d :\; N(\d )\neq C(\d )\neq D(\d ), N(\d )+D(\d )=F(\d )\neq A\}$, a set of
{\em weakly nil-semi-simple} derivations.
\end{enumerate}
For each element $a\in A$, one can attach, so-called, the {\em
inner derivation} $\ad \, a$ of the algebra $A$ by the rule $\ad
\, a(b):=[a,b]=ab-ba$, and then four subalgebras of $A$: the {\em
centralizer} $C(a):=C(\ad \, a)=\{ b\in A\, | \, ab=ba\}$ of the
element $a$ in $A$, $N(a):=N(\ad \, a)$,  $D(a):=D(\ad \, a)$, and
$F(a):=F(\ad \, a)$. The algebra $A$ is a Lie algebra where the
bracket is the commutator of elements. The map $\ad :A\ra
\Der_K(A)$, $a\ra \ad \, a$, is a Lie algebra homomorphism with
kernel $Z(A)$, the {\em centre} of the algebra $A$.

The set $A\backslash Z(A)$ is a disjoint union of its eleven
subsets
$$ A\backslash Z(A) =\bigcup_{i=1}^{11}\, \D_i(A)\;\; {\rm where}\;\;
\D_i (A):=\{ a\in A\backslash Z(A)\,\, |\,\, \ad \, a\in
\widetilde{\Delta}_i(A)\}.  $$

Elements of the sets $\D_1(A), \ldots , \D_{11}(A)$ are called
respectively {\em strongly nilpotent}, $\ldots $,  {\em weakly
nil-semi-simple}. The sets $\D_i(A)$ are invariant under the
action of the group $\Aut_K(A)$ of all $K$-algebra automorphisms
of $A$.

{\sc Definition}. We say that an  algebra $A$ admits the {\bf
Dixmier partition} if $$ A\backslash Z (A)=\D_1 (A)\cup \D_2
(A)\cup \D_3 (A)\cup \D_4(A)\cup \D_5 (A).$$ Equivalently, an
algebra $A$ admits the Dixmier partition if and only if, for each
$a\in A\backslash Z(A)$, either $F(a)=N(a)$ or $F(a)=D(a)$ (the
{\em nilpotent-semi-simple alternative}).

 Dixmier \cite{Dix}
proved that the first Weyl algebra $A_1$ satisfies this property.
 One can modify arguments of Dixmier and prove that {\em each
 non-commutative algebra $A$ of Gelfand-Kirillov dimension $<3$
 admits the Dixmier partition provided $\bA $ is a domain} (Corollary \ref{2NorD}).

{\bf Dixmier's problem 6}, \cite{Dix}: {\em If $a\in \D_5(A_1)$
then $K[a]\backslash K\subseteq \D_5(A_1) $ ?}

The type (properties) of an element $a\in A_1$   is determined by
the type of the corresponding inner derivation $\ad \, a$. So,
Dixmier's problem 6 is a question about properties of inner
derivations. Dixmier proved in his paper \cite{Dix}, 10.3, that,
{\em for $i=1,2$ and any non-scalar  polynomial} $f(t)$, $a\in
\D_i(A_1)$ $ \Leftrightarrow $ $f(a)\in \D_i(A_1)$. So,   for an
arbitrary algebra $A$, Dixmier's problem  6 is essentially a {\em
question} of
\begin{itemize}
\item {\em whether an inner derivation of the type $\ad \, f(a)$
has a non-zero eigenvalue (where $a\in A$, $f(t)\in K[t]$, $\deg\,
f(t)>1$). }
\end{itemize}
The answer is {\em negative} for many classes of algebras (rings
of differential operators on  smooth irreducible algebraic
varieties, universal enveloping algebras  of Lie algebras) as
follows from the next result (Section  \ref{noneadfa}, see Theorem
\ref{Beigen} for a more general result).

\begin{theorem}\label{Beigvfa=0}
 Let $B=\cup_{i\geq
0}\, B_i$ be a filtered algebra over a field $K$ of characteristic
zero such that the associated graded algebra $\gr\, B$ is a
commutative domain, and let $a\in B$.
 Then, for an arbitrary polynomial $f(t)\in
K[t]$ of degree $>1$, $0$ is the only eigenvalue of the inner
derivation $\ad\, f(a)$ of the algebra $B$ (that is, if $\d :=\ad
\, f(a)\neq 0$ and $\d (b)=\l b$ for some $\l \in K$ and $0\neq
b\in B$ then $\l =0$).
\end{theorem}

As an immediate consequence of Theorem \ref{Beigvfa=0} we have the
following result (among others of the same sort, Section
\ref{noneadfa}) about {\em nonexistence of nonzero eigenvalue  for
inner derivations of the type $\ad \, f(a)$}.

\begin{theorem}\label{CDXnonex}
Let $\CD (X)$ be the ring of differential operators on a smooth
irreducible affine algebraic variety $X$ over a field $K$ of
characteristic zero, and let $a\in \CD (X)$. Then, for an
arbitrary polynomial $f(t)\in K[t]$ of degree $>1$, $0$ is the
only eigenvalue of the inner derivation $\ad\, f(a)$ of $\CD (X)$.
\end{theorem}

{\it Proof} (see Section \ref{noneadfa} for detail).
 By the very definition of the ring of differential operators $\CD (X)$, the algebra $\CD
 (X)$ has the order filtration $ \CD (X)=\cup_{i\geq 0} \CD (X)_i$
  (by the total degree of derivations)
 such  that the associated graded algebra $\gr \, \CD (X)$ is a
 commutative domain, \cite{MR}, 15.4.7 and 15.5.6. So, the result
 follows from Theorem \ref{Beigvfa=0}. $\Box $

The same result is true $(i)$  for all prime factor algebras of
the universal enveloping algebra $U(\CG )$ where $\CG $ is a
completely solvable, algebraic Lie algebra (Corollary
\ref{Gcomplsolv}); or $\CG $ is a nilpotent Lie algebra (Corollary
\ref{noeignilpr}); $(ii)$ for all primitive factor algebras of the
universal enveloping algebra $U(\CG )$ where $\CG $ is a solvable
Lie algebra (Corollary \ref{Ggencs}); $(iii)$ for algebra $\CA (V,
\d  , \G )$ from \cite{MR}, Ch. 14, Section 8.

Now we are ready to  give a short proof to Dixmier's problem 6.

{\bf Proof (Dixmier's problem 6)}. We may assume that $K=\bK $
(since $\bA_1$  is a domain). Recall that the Weyl algebra $A_1$
satisfies the Dixmier partition \cite{Dix}.  The Weyl algebra
$A_1$ has a filtration by the total degree of the canonical
generators:
$$A_1=\bigcup_{n\geq 0}\, A_{1,n}\;\; {\rm where}\;\;
A_{1,n}:=\sum_{i,j\geq 0}\{  Kx^i\der^j\, | \, i+j\leq n\}.$$ The
associated graded algebra $\gr \, A_1:=\oplus_{n\geq 0} \, A_{1,
n}/A_{1, n-1}$ is a polynomial algebra $K[ \overline{x},
\overline{\der}]$ in two variables $ \overline{x}:=x+A_{1,0}$ and
$ \overline{\der}:=\der+A_{1,0}$ where $A_{1,0}=K$ and
$A_{1,-1}:=0$.

Let $a\in \D_5(A_1)$ and $f(t)\in K[t]$, $\deg \, f(t)>1$. By
Theorem \ref{Beigvfa=0}, $f(a)\not\in  \D_{3,4}(A_1)$. By the
result of Dixmier \cite{Dix}, 10.3 (mentioned above), $f(a)\not\in
\D_{1,2}(A_1)$. Hence $f(a)\in \D_5(A_1)$. $\Box $

Dixmier's question $6$ is the most difficult part in an answer to
the question: {\em suppose that the type of an element $a$ is
known. What is the type of the element $f(a)$ where $f(t)\in
K[t]$, $\deg_t \, f(t)>1$?} For the Weyl algebra $A_1$ the answer
is
\begin{enumerate}
\item {\em Let $i=1,2,5$. If $a\in \D_i(A_1)$ then
$f(a)\in \D_i(A_1)$.
\item Let $i=3,4$. If $a\in \D_i(A_1)$ then
$f(a)\in \D_5(A_1)$.}
\end{enumerate}

{\it Proof}. For $i=1,2$, this is the result of Dixmier
\cite{Dix}, 10.3. For $i=5$, this is Dixmier's question $6$. For
$i=3,4$, $f(a)\not\in \D_{3,4}(A_1)$ (see the proof to Dixmier's
problem $6$), and $f(a)\not\in \D_{1,2}(A_1)$, by \cite{Dix},
10.3, and so $f(a)\in \D_5(A_1)$. $\Box $

Exactly the same result holds for the ring of differential
operators $\CD (X)$ on a smooth irreducible algebraic curve $X$
(Corollary \ref{CDXDP}), or even in more general situation
(Theorem \ref{DP6ccc}).

The elements $a$ and $f(a)$ {\em commute}. One can ask a similar
question but for commuting elements.

{\bf Question}. {\em Suppose that the type of an element $a$ is
known. What are the types of elements commuting with $a$?}

For many classes of algebras an answer is given by the following
definition (see Theorem \ref{G3Apr} and the examples that follow).

{\sc Definition}. Let $A$ be a noncommutative  algebra that admits
the Dixmier partition. We say that $A$ satisfies the {\bf
homogeneous centralizer condition} (the {\em hcc}, for short) if,
for each element $a\in A \backslash Z(A)$, all the elements of the
set $C(a)\backslash Z(A)$ have exactly one of the following types:
\begin{enumerate}
\item strongly nilpotent,
\item  weakly nilpotent,
\item  generic except for a subset $K^*a+Z(A)$  of strongly
semi-simple  elements for some $a\in \D_3(A)$ where
$K^*:=K\backslash \{ 0\}$,
\item generic except for a subset $K^*a+Z(A)$  of weakly
semi-simple  elements for some $a\in \D_4(A)$,
\item generic.
\end{enumerate}

The next  definition is inspired by the result of Flanders (see
the review of Ore of the paper \cite{Ami} on MathSciNet) and
Amitsur \cite{Ami}: {\em Let $a$ be a non-scalar element of the
Weyl algebra $A_1$. Then its centralizer $C(a)$ is a commutative
subalgebra of $A_1$, and a free, finitely generated
$K[a]$-module.}

{\sc Definition}.  A noncommutative algebra $A$ satisfies the {\bf
commutative centralizer condition} (the $ccc$, for short) if
 the  centralizer of each element of the
set $A\backslash Z(A)$ is a {\em commutative} algebra.

When this paper was written D. Jordan and T. Lenagan pointed out
on a paper of K. R. Goodearl on centralizers \cite{goodearlcen}
where the reader can find numerous examples of algebras that
satisfy the commutative centralizer condition.

Suppose that an  algebra $A$  satisfies the commutative
centralizer condition. Then, for each element $a\in A \backslash
Z(A)$, its centralizer $C(a)$ is a unique maximal commutative
subalgebra of $A$ that contains the element $a$, and each maximal
 commutative subalgebra $C$ of $A$ coincides with the centralizer of
every element of $C\backslash Z(A)$ (Corollary
\ref{PAmaxcomsubr}). So, for an
 algebra that satisfies  the {\em ccc} the concepts of
centralizer and of maximal isotropic subalgebra {\em coincide}.
For $a, b\in A\backslash Z(A)$, either $C(a)=C(b)$ or otherwise
$C(a)\cap C(b)=Z(A)$ (Corollary \ref{PACf=Cg}),  and so every
algebra that satisfies the {\em ccc} is a union of distinct
centralizers ($=$ maximal commutative subalgebras) that meet at
$Z(A)$. If, in addition, the algebra $A$ admits the Dixmier
partition and satisfies the {\em hcc}, then this union $A=\cup \,
C(a)$ has a transparent (homogeneous) structure. One of the goals
of this paper is to prove that many popular algebras admit the
Dixmier partition and satisfy both the homogeneous centralizer
condition and the commutative centralizer condition.

The following statement is one of the  main results of the  paper
(Section \ref{ccc}).

\begin{theorem}\label{G3Apr}
Let a $\bK$-algebra $\G $ be a division ring that satisfies the
commutative centralizer condition, and let  $A$ be a
non-commutative  $K$-algebra of Gelfand-Kirillov dimension $<3$
over $K$  such that $\bA$ is (isomorphic to) a $\bK $-subalgebra
of  $\G $.
 Then $A$ admits the Dixmier partition and satisfies both the
 homogeneous centralizer condition and the commutative centralizer
 condition.
\end{theorem}

Let $\CD (X)$ be the ring of differential operators on a smooth
irreducible algebraic curve $X$ over $K=\bK $. Then the division
ring $ Q (\CD (X))$ of the algebra $\CD (X)$ satisfies the {\em
ccc} (Corollary \ref{CDGccc}). Let $\CL =L[[X,X^{-1};\s ]]$ be the
{\em skew Laurent series algebra} with coefficients from a field
$L$ where an automorphism  ${\rm id}_L\neq \s \in \Aut (L)$
satisfies the following condition: for each $l\in L$ either $\s
(l)=l$ or $\s^i (l) \neq l$ for all $i\geq 1$. Let $\CR
=L((t^{-1};\d ))$ be the {\em formal pseudo-differential operator
ring} with $0\neq \d \in \Der (L)$. Then the division rings $\CL $
and $\CR $ satisfy the
 {\em ccc} (Theorem \ref{CLccc} and Lemma \ref{DSccr}). The last
 fact is due to K. R. Goodearl \cite{goodearlcen}.  Many well-known
algebras of Gelfand-Kirillov dimension $<3$ are subalgebras of
either $\CL $ or $\CR $, and so satisfy the {\em ccc} as well.

Using the result above we prove {\em that the following algebras
admit the Dixmier partition and satisfy  the commutative
centralizer condition and the homogeneous centralizer condition:}
\begin{enumerate}
\item The first Weyl  algebra $A_1$. Any noncommutative
subalgebra $A$ of the division ring $Q(A_1)$ of the Weyl algebra
$A_1$ with Gelfand-Kirillov dimension $\GK (A)<3$. In particular,
the {\em noncommutative deformations of type-A Kleinian
singularities} are of this type (Corollary \ref{QA1ch}).
\item The {\em quantum plane} $\L =K\langle x,y\,| \, xy=\l
yx\rangle $, $\l \in K^*$ is not an $i^{th}$ root of $1$ for all
$i\geq 1$. Any noncommutative subalgebra $A$ with Gelfand-Kirillov
dimension $\GK (A)<3$ of the division ring $Q(\L )$ of the quantum
plane $\L $ (Corollary \ref{QLch}). The fact that $ Q(\L )$
satisfies the {\em ccc} is due to V. A. Artamonov and P. M. Cohn
\cite{Art-Cohnqpl}. The {\em quantum Weyl } algebra $A_1(\l
)=K\langle x, y \, | \, yx-\l xy =1 \rangle$ is a subalgebra of
$Q(\L )$, hence it satisfies the {\em ccc}, this was proved by V.
Mazorchuk \cite{Mazccc}.

\item $(K=\bK)$ The ring $\CD (X)$ of
differential operators on a smooth  irreducible algebraic curve
$X$. Any noncommutative subalgebra $A$ with Gelfand-Kirillov
dimension $\GK (A)<3$ of the division ring $Q(\CD (X) )$ of
 the algebra $\CD (X)$ (Corollary \ref{AsubDX}).

\item The universal enveloping algebra $Usl(2)$ of the Lie algebra
$sl(2)$ (note that $\GK (Usl(2))=3)$. Any noncommutative
subalgebra $A$ of the division ring $Q(Usl(2))$ of the universal
enveloping algebra $Usl(2)$ with Gelfand-Kirillov dimension $\GK
(A)<4$  and such that $A\cap  K(C)\neq K$ where $C$ is the {\em
Casimir} element of the algebra $Usl(2)$ (Proposition
\ref{Usl2p3}).

\item The quantum $U_qsl(2)$ (note that $\GK (U_qsl(2))=3)$. Any noncommutative
subalgebra $A$ of the division ring $Q(U_qsl(2))$ of the  algebra
$U_qsl(2)$ with Gelfand-Kirillov dimension $\GK (A)<4$  and such
that $A\cap  K(C)\neq K$ where $C\in U_qsl(2)$ is the (quantum)
{\em Casimir} element (Proposition \ref{Uqsl2p}).
\end{enumerate}



\section{Nonexistence of Nonzero Eigenvalue for Inner Derivations
of the Type $\ad\, f(a)$, $f(t)\in K[t]$, $\deg\,
f(t)>1$}\label{noneadfa}

Let $K$ be a field of characteristic zero. In this Section, let
$B=\cup_{i\geq 0}\, B_i$ be a filtered $K$-algebra such that
 the associated graded algebra $\gr\,
B=\oplus_{i\geq 0}\, B_i/B_{i-1}$ is a {\em domain}, and
$$ [B_i,B_j]\subseteq B_{i+j-1}\;\; {\rm for \; all}\;\; i,j\geq
0, \;\; i+j>0.$$ In this case, we say that the algebra $B$ has a
{\em filtration of nilpotent type}, and $\{ B_i\}$ is the {\em
filtration of nilpotent type}. If, in addition, all the elements
of the set $B_0\backslash Z(B)$ are strongly nilpotent then we say
that the algebra $B$ has a filtration of {\em strongly nilpotent
type}, and $\{ B_i\}$ is a filtration of {\em strongly nilpotent
type}. Note that, in general, the algebra $\gr \, B$ is not
necessarily commutative. Although in applications this is usually
the case. If $\gr\, B$ is a commutative domain then the filtration
$\{ B_i\}$ is of strongly nilpotent type.

 For a nonzero element
$a\in B$ the unique natural number $n$ such that $a\in
B_n\backslash B_{n-1}$ is called the {\em degree} of the element
$a$ (with respect to the filtration) denoted $\deg \, a$ $(\deg\,
0:=-\infty )$. The algebra $\gr \, B$ is a  domain, thus, for any
$a, b\in B$,
\begin{eqnarray*}
\deg (ab) & = &\deg\, a+\deg\,b ,\\
\deg(a+b)&\leq &\max \{ \deg\, a, \deg\,b\},\\
\deg\, [a,b] & \leq &\deg\, a+\deg\,b -1, \;\; {\rm if}\;\; \deg\, a+\deg\, b>0.\\
\end{eqnarray*}
The first property guarantees that $B$ is a domain since $\gr\, B$
is so.

For a given element $a\in B$, the following three linear maps from
$B$ to $B$ commute: the left multiplication $L_a$ by the element
$a$, the right multiplication $R_a$ by the element $a$, and the
inner derivation $\ad \, a=L_a-R_a$ of $B$. So, for each $n\geq
1$,
$$ (L_a)^n=(\ad \, a +R_a)^n=\sum_{i=0}^n \, {n\choose i}(\ad \,
a)^i\, (R_a)^{n-i}. $$ Thus, for any $b\in B$,
\begin{eqnarray*}\label{anb=Lanb}
a^nb & = &(L_a)^nb=(\ad \, a +R_a)^nb=\sum_{i=0}^n \, {n\choose i}
(\ad \, a)^i(b)\, a^{n-i}\\
& = & \sum_{i=0}^n\, n(n-1)\cdots (n-i+1)  \, {(\ad \,
a)^i(b)\over i!}\, a^{n-i}.
\end{eqnarray*}

For a polynomial $f(t)\in K[t]$ we denote by $f^{(i)}(t)$ its
$i^{th}$ derivative ${d^if\over dt^i}$ with respect to $t$.  If
$\deg\, a >0$ and $\deg_t (f)>0$  then it follows from the
identity above that
\begin{equation}\label{comfab}
[f(a),b]=\sum_{i\geq 1}\, {(\ad \, a)^i(b)\over i!}\, f^{(i)}(a)=
(\ad \, a)(b)\, f'(a)+\cdots
\end{equation}
where the dots denote the terms of smaller degree. Moreover,
$$ \deg \, T_1> \deg \, T_2 >\cdots , \;\;
{\rm where}\;\;  T_i:={(\ad \, a)^i(b)\over i!}\, f^{(i)}(a).$$ In
fact, if $T_i\neq 0$ and $i>1$ then (using $\deg \, a>0$)
\begin{eqnarray*}
\deg \, T_i &=&\deg \, (\ad \, a)^i(b)+\deg \,f^{(i)}(a)\\
&=&\deg \, (\ad \, a ) (\ad \, a)^{i-1}(b)+\deg \,f^{(i-1)}(a)
-\deg
\, a\\
&\leq & \deg \, (\ad \, a)^{i-1}(b)+\deg (a) -1+\deg
\,f^{(i-1)}(a) -\deg \, a\\
&=&\deg \, T_{i-1}-1 <\deg \, T_{i-1}.
\end{eqnarray*}

\begin{theorem}\label{Beigen}
Let $B=\cup_{i\geq 0}\, B_i$ be a filtered $K$-algebra such that
 the associated graded algebra $\gr\,
B$ is a domain, and $ [B_i,B_j]\subseteq B_{i+j-1}$  for  all
$i,j\geq 0$, $i+j>0$. Let $a$ be an element of the algebra $B$
such that, for some polynomial $f(t)\in K[t]$ of degree $>1$, the
inner derivation $\ad \, f(a)$ of the algebra $B$ has a nonzero
eigenvalue. Then $a_0\in B_0$.
\end{theorem}

{\it Proof}. By the assumption, $[f(a),b]=\l b$ for some $0\neq \l
\in K$ and $0\neq b\in B$. Without loss of generality we may
assume that $\l =1$. Suppose that the element $a$ does not belong
to $B_0$ thus $\deg \, a>0$. We seek a contradiction.
 By (\ref{comfab}),
$$ b=[f(a),b]=\sum_{i\geq 1}\, {\d^i(b)\over i!}\, f^{(i)}(a),
\;\; {\rm where}\;\; \d :=\ad \, a.$$ Observe that $\d
(f^{(i)}(a))=[a,f^{(i)}(a)]=0$ for all $i$, thus, if we substitute
the sum for $b$ above in the right-hand side of the equality above
we obtain
$$ b=\sum_{i_1, i_2\geq 1}\, {\d^{i_1+i_2}(b) \over
i_1!i_2!}\,f^{(i_1)}(a)f^{(i_2)}(a).$$ Repeating this substitution
$s$ times we have

\begin{equation}\label{b=bigs}
 b=\sum_{i_1,\ldots , i_s \geq 1}\,
{\d^{i_1+\cdots +i_s}(b) \over i_1!\cdots
i_s!}\,f^{(i_1)}(a)\cdots f^{(i_s)}(a).
\end{equation}
Now, we can find the leading term of the sum above,
\begin{equation}\label{b=lts}
 b=\d^s(b)[f'(a)]^s+\cdots .
\end{equation}
The proof of (\ref{b=lts}) is essentially the same as in the case
when $s=1$, and this result follows immediately from the next
fact: {\em let} $T_{(i_1,\ldots , i_s)}:= \d^{i_1+\cdots
+i_s}(b)\,f^{(i_1)}(a)\cdots f^{(i_s)}(a)${\em ; then}
$$ \deg\, T_{(i_1+1,i_2,\ldots , i_s)}<\deg\, T_{(i_1,i_2,\ldots ,
i_s)}, \;\;  {\it provided} \;\; T_{(i_1+1,i_2,\ldots , i_s)}\neq
0;$$ since
\begin{eqnarray*}
\deg\, T_{(i_1+1,i_2,\ldots , i_s)} & = & \deg\, \d \d^{i_1+\cdots
+ i_s}(b)+\deg\, f^{(i_1+1)}(a)+\sum_{\nu \geq 2}\, \deg\,
f^{(i_\nu )}(a)\\
& \leq & \deg\, \d^{i_1+\cdots +i_s}(b)+\deg (a) -1+ \deg\,
f^{(i_1)}(a)-\deg (a) +\sum_{\nu \geq 2}\, \deg\, f^{(i_\nu )}(a)\\
& = &T_{(i_1,i_2,\ldots , i_s)}-1 <\deg\, T_{(i_1,i_2,\ldots ,
i_s)}.
\end{eqnarray*}
Clearly, $\d^n(b)\neq 0$ for all $n\geq 1$ (since otherwise, by
(\ref{b=bigs}), we would have $b=0$, a contradiction). Now, by
(\ref{b=lts}) and since $\deg_t(f'(t))\geq 1$, $\deg (a)\geq 1$,
we obtain 
\begin{equation}\label{degbgs}
\deg\, b=\deg\, \d^s(b)+s\deg\, f'(a)\geq s\, \deg_t(f'(t))\cdot
\deg (a)\geq s, \;\; {\rm for \;\; all}\;\; s\geq 1,
\end{equation}
 which is impossible since $\deg \, b<\infty $. This contradiction finishes
 the proof of the theorem.
$\Box $

{\bf Proof of Theorem \ref{Beigvfa=0}}.

The associated graded algebra $\gr\, B$ is commutative which
implies that the elements $B_0\backslash Z(B)$ are strongly
nilpotent $([B_i,B_j]\subseteq B_{i+j-1}$ for all $i,j\geq 0$, and
so $(\ad\, B_0)^{i+1}B_i=0)$. Suppose that, for some element $a\in
B$ and some polynomial $f(t)$ of degree $>1$, the inner derivation
$\ad \, f(a)$ of the algebra $B$ has a nonzero eigenvalue. Then
$f(a)\not\in Z(B)$ and,  by Theorem  \ref{Beigen}, $a\in B_0$ thus
$f(a)\in B_0\backslash Z(B)$ which contradicts  to the fact that
each element of $B_0\backslash Z(B)$ is  strongly nilpotent. $\Box
$

 Let $C$ be a {\em finitely
generated} algebra over the field $K$. To each finite set of
$K$-algebra generators $x_1, \ldots , x_n$ one can attach,
so-called, the {\em standard} filtration of the algebra $C$ (by
the  {\em total} degree of generators):
$$ C_0:=K\, \subseteq \, C_1:=K+\sum_{j=1}^n\, Kx_j\, \subseteq \cdots
\subseteq \, C_i:=C_1^i\, \subseteq \cdots .$$

{\sc Definition}. A finitely generated algebra $C$  is called an
{\em almost commutative} algebra if there exists  a standard
filtration in $C$ such that the associated graded algebra $\gr\,
C=\oplus_{i\geq i}\, C_i/C_{i-1}$ is a {\em commutative} algebra.

{\em Each} factor algebra of the universal enveloping algebra
$U(\CG )$ of a finite dimensional Lie algebra $\CG $ is an almost
commutative algebra and vice versa (\cite{MR}, 8.4.3, or
\cite{KL}).

 {\it Example 1}  ({\bf Almost Commutative Algebras}). Any almost commutative algebra $C$
such that the associated graded algebra $\gr\, C$ is a domain
satisfies the conditions of Theorem \ref{Beigvfa=0} with $C_0=K$.

{\it Example 1.} $(i)$. The natural filtration of the universal
enveloping algebra $U(\CG )$ is a filtration of strongly nilpotent
type since $U(\CG )_0=K$ and the associated graded algebra $\gr \,
U(\CG )$ is a polynomial algebra in $\dim_K \, \CG $
indeterminates.

{\it Example 1.} $(ii)$  ({\bf the Weyl Algebras}). The {\em Weyl
algebra} $A_n$ is an associative algebra generated over the field
$K$ by $2n$ elements
 $X_1, \ldots , X_n, \der_1, \ldots , \der_n$ subject to the
 defining relations:
 \begin{equation*}\label{defrelweyl-1}
\der_iX_i-X_i\der_i=1, \;\; i=1, \ldots , n,
\end{equation*}
\begin{equation*}\label{defrelweyl-2}
[X_i, X_j]=[\der_i, \der_j]=[\der_i, X_j]=0, \;\; {\rm for \;
all}\;\; i\neq j.
\end{equation*}
The Weyl algebra $A_n$ is canonically isomorphic to the ring of
differential operators
$$A_n\simeq K[X_1, \ldots , X_n, {\der \over \der X_1}, \ldots , {\der \over \der
X_n}],\;\; X_i\lra X_i, \;\; \der_i\lra {\der \over \der X_i},$$
with polynomial coefficients $K[X_1, \ldots , X_n]$. The field $K$
has characteristic zero, so the centre of the Weyl algebra $A_n$
is $K$. The standard filtration of the Weyl algebra $A_n$
associated with the canonical generators is called the {\em
Bernstein} filtration. The associated graded algebra of the Weyl
algebra $A_n$ with respect to this filtration is a polynomial
algebra in $2n$ indeterminates:
$$ K[x_1, \ldots , x_n, y_1, \ldots , y_n], \;\; x_i:=X_i+K\in A_{n,1}/A_{n,0}, \;\;
y_i:=\der_i+K\in A_{n,1}/A_{n,0}.$$ Thus the Bernstein filtration
of
 the Weyl algebra $A_n$ is a filtration of strongly nilpotent
 type.
\begin{corollary}\label{Anfaev=0}
Let $a$ be a non-scalar element of the Weyl algebra $A_n$ over a
field $K$ of characteristic zero. Then, for an arbitrary
polynomial $f(t)\in K[t]$ of degree $>1$, $0$ is the only
eigenvalue of the inner derivation $\ad\, f(a)$ of $A_n$. $\Box $
\end{corollary}

An {\em arbitrary} $K$-derivation of the Weyl algebra $A_n$ over a
field of characteristic zero is an {\em inner} derivation
\cite{Dix3}, so  Corollary \ref{Anfaev=0} is a statement about
derivations of the Weyl algebra. In particular, if an (inner)
derivation $\d =\ad\, b$ of the Weyl algebra $A_n$ has a nonzero
eigenvalue then the element $b\in A_n$ is {\em not} a polynomial
$f(a)$ of any element $a\in A_n$.

%
%

{\it Example 2}  ({\bf Rings of Differential Operators on Smooth
Algebraic Varieties}). Let $R$ be a commutative algebra over  a
field $K$ of characteristic zero. The {\em ring of differential
operators} $\CD (R)=\cup_{i\geq 0}\, \CD (R)_i$ with coefficients
from the algebra $R$ is a positively filtered subalgebra of the
algebra $\End_K(R)$ of $K$-linear maps in $R$ defined by the rule
$$ \CD(R)_0:=\{ f\in \End_K(R)\, | \, fr-rf=0\;\; {\rm for \;\; all}\;\; r\in R\}=
\End_R(R)=R $$ and
$$\CD(R)_i:=\{ f\in \End_K(R)\, | \, fr-rf\in \CD(R)_{i-1}\;\;
{\rm for \;\; all}\;\; r\in R\}.$$ For a smooth irreducible affine
variety $X$ over the field $K$ its coordinate algebra $R=\OO (X)$
 is a commutative affine regular domain. In this case, the ring of
differential operators $\CD (X):=\CD (\OO (X))$ on $X$ is
generated, as a $K$-subalgebra of $\End_K(\OO (X))$, by the
coordinate  algebra $\OO (X)$ and the $\OO (X)$-module $\Der_K(\OO
(X))$ of $K$-derivations of the algebra $\OO (X)$ (\cite{MR},
15.5.6). The associated graded algebra $\gr\, \CD
(X)=\oplus_{i\geq 0} \, \CD (X)_i/\CD (X)_{i-1}$ is a commutative
(affine) domain, \cite{MR}, 15.4.7.

\begin{corollary}\label{noeignilpr}
Let $\CG$ be a nilpotent finite dimensional Lie algebra over a
field $K$ of characteristic zero, $A$ be a noncommutative prime
factor ring of its universal enveloping algebra $U(\CG )$, and
$a\in A$ be an arbitrary element. Then, for an arbitrary
polynomial $f(t)\in K[t]$ of degree $>1$, $0$ is the only
eigenvalue of the inner derivation $\ad\, f(a)$ of $A$.
\end{corollary}

{\it Proof}. The algebra $A$ is a domain such that the
localization $S^{-1}A$ of the algebra $A$ at $S:=Z(A)\backslash \{
0\} $ is the Weyl algebra $A_n(Q)$ with coefficients from the
field $Q:=S^{-1}Z(A)$ of quotients of the centre $Z(A)$ of $A$
(\cite{MR}, 14.6.9.(i)). Now, the result follows from Corollary
\ref{Anfaev=0}. $\Box $

We denote by $A_1'$ the localization of the Weyl algebra $A_1$ at
the powers of the element $x$. Since $A_1=\CD (K[x])$, we have
$A_1'=\CD (K[x,x^{-1}])$. For a natural number $m$, let
$A_m':=A_1'\t \cdots \t A_1'$ , $m$ times.

\begin{corollary}\label{CAnAm}
Let $A=C\t A_n\t A_m'$ where $C$ is a commutative algebra and
$n,m\geq 1$, and let $a\in A$. Then, for an arbitrary polynomial
$f(t)\in K[t]$ of degree $>1$, $0$ is the only eigenvalue of the
inner derivation $\ad\, f(a)$ of $A$.
\end{corollary}

{\it Proof}. Note that the algebra $A$ is a domain. On the one
hand, the algebra $A$ is the Weyl algebra $A_n(\L )$ with
coefficients from the algebra $\L :=C\t A_m'$. The Bernstein
filtration $A_n(\L )=\cup_{i\geq 0}A_n(\L )_i$ for the Weyl
algebra $A_n(\L )$ is a filtration of nilpotent type since $ \gr\,
A_n(\L )$ is isomorphic to the polynomial algebra $\L [x_1, \ldots
, x_n , y_1, \ldots , y_n]$  in $2n$ indeterminates with
coefficients from the algebra $\L $.

On the other hand, the algebra $A=\L'\t A_m'=\L'\t \CD (Q)$ where
$\L':=C\t A_n$ and $Q:=K[z_1, z_1^{-1}, \ldots , z_m,z_m^{-1}]$ is
the Laurent polynomial algebra. The order filtration $A_m'=\CD
(Q)=\cup_{i\geq 0}\CD (Q)_i$ of the algebra  $A_m'$ can be
extended to the filtration $A=\cup_{i\geq 0}B_i$, $B_i:=\L' \t \CD
(Q)_i$, for the algebra $A$. The filtration  $\{ B_i\}$ is of
nilpotent type since $$\gr \, A\simeq \L'\t \gr\, \CD (Q)\simeq
\L'\t K[z_1, z_1^{-1}, \ldots , z_m,z_m^{-1}, y_1, \ldots ,
y_m].$$ Note that $A_n(\L )_0=\L =C\t A_m'$, $B_0=\L'\t \CD
(Q)_0=C\t A_n\t Q$, and $I:=A_n(\L )_0\cap  B_0=C\t Q$. Each
element of the set $I\backslash Z(A)$ is a strongly nilpotent
element of the algebra $A$.

Suppose that, for some element $a\in A$ and some polynomial $f(t)$
of degree $>1$, the inner derivation $\ad \, f(a)$ of the algebra
$A$ has a nonzero eigenvalue. Then $f(a)\not\in Z(A)$ and, by
Theorem \ref{Beigen}, $a\in I$, a contradiction. $\Box $

Recall that an {\em algebraic} Lie algebra $\CG $ is the Lie
algebra of a (linear) affine algebraic  group.  The reader is
referred to the books \cite{BGR} and \cite{MR}  for detail.

\begin{corollary}\label{Gcomplsolv}
Let $\CG$ be a completely solvable, algebraic Lie algebra over a
field $K$ of characteristic zero, $A$ be a prime factor algebra of
the universal enveloping algebra $U(\CG )$ and $a$ be an arbitrary
element of the algebra $A$. Then, for an arbitrary polynomial
$f(t)\in K[t]$ of degree $>1$, $0$ is the only eigenvalue of the
inner derivation $\ad\, f(a)$ of $A$.
\end{corollary}

{\it Proof}. By \cite{MR}, 14.2.11, $A$ is a domain. By \cite{MR},
14.7.6, the localization $A_e$ of the algebra $A$ at the powers of
the particularly chosen element $e\in A$ is the  algebra from
Corollary \ref{CAnAm}, $A_e\simeq Z(A_e)\t A_n\t A_m'$. So, the
result follows from Corollary \ref{CAnAm}. $\Box $

Let $A$ be an algebra. A left $A$-module $M$ is called {\em
faithful} if $aM=0$  implies $a=0$ for any $a\in A$. If an algebra
$A$ has a left faithful simple module then $A$ is called {\em
(left) primitive}.

\begin{corollary}\label{Ggencs}
Let $\CG$ be a completely solvable  Lie algebra over a field $K$
of characteristic zero, $A$ be a prime factor algebra of the
universal enveloping algebra $U(\CG )$,  and $a\in A$. Then, for
an arbitrary polynomial $f(t)\in K[t]$ of degree $>1$, $0$ is the
only eigenvalue of the inner derivation $\ad\, f(a)$ of $A$. $\Box
$
\end{corollary}

{\it Proof}. By \cite{MR}, 14.2.11, $A$ is a domain. Combining
 \cite{MR}, 14.9.16 and 14.8.11, $A$ is a subalgebra of the
 algebra from  Corollary \ref{CAnAm} with $C$ a polynomial
 algebra. Now, the result  follows from Corollary \ref{CAnAm}. $\Box $

{\it Remark}. The same result holds for the algebra $\CA (V, \d ,
\G )$ from \cite{MR}, Ch. 14, Section 8 (since, by \cite{MR},
14.8.11, $\CA (V, \d , \G )$ is a subalgebra of the algebra from
Corollary \ref{CAnAm} with $C$ is a polynomial algebra).

\begin{corollary}\label{cor11}
Let the algebra $B$ be as in Theorem \ref{Beigvfa=0}, and let
$a\in B\backslash Z(B)$. If $f(a)\in Z(B)$ for some polynomial
$f(t)\in K[t]$ of degree $>1$ then $0$ is the only eigenvalue for
the inner derivation $\ad (a)$, and $a\in B_0$. If, in addition,
$B_0\subseteq Z(B)$ then $f(x)\in Z(B)$, $x\in B$, imply $x\in
Z(B)$.
\end{corollary}

{\it Proof}. $\l b=[a,b]=[f(a)+a, b]$, $0\neq b\in B$, $\l \in K$,
imply $\l =0$, by Theorem \ref{Beigvfa=0}.

Suppose that $a\not\in B_0$. Then, by (\ref{comfab}), $0=[f(a),
b]=[a,b]f'(a)+\cdots $ for all $b\in B$, hence $[a,b]=0$ for all
$b\in B$, that is $a\in Z(B)$, a contradiction. Therefore, $a\in
B_0$. Then the last statement is obvious. $\Box $


\section{Noncommutative Domains of Gelfand-Kirillov $<3$
Admit the Dixmier Partition}\label{Dom3DP}

In this section, if it is not stated otherwise, $K$ is {\em an
algebraically closed field of characteristic zero}, and a
$K$-algebra $A$ is a {\em noncommutative domain}.

The aim  of this section is to prove that the algebra $A$ admits
the Dixmier partition provided that its Gelfand-Kirillov dimension
is less than $3$ (Corollary \ref{NorD}).

For  each non-central  element $a\in A \backslash Z(A)$, we have
the corresponding nonzero inner derivation $\ad \, a$ of the
algebra $A$. Its torsion algebra
$$ F(a)=\bigoplus_{\l \in K}F(a, \l ), \;\;  F(a, \l )=\cup_{n\geq
0}F(a,\l  , n),$$ where $ F(a,\l  , n):=\ker (\ad \, a
-\l)^{n+1}=\{ x\in A\, | \, (\ad \, a -\l)^{n+1}x=0\}$. For any
$x,y\in A$, $\l , \mu \in K$, and $n\geq 1$, the identity
\begin{equation}\label{dlmn}
(\ad\, a -(\l +\mu ) )^{n}(xy)=\sum_{i=0}^n\, {n\choose i}(\ad \,
a -\l )^i(x)\cdot (\ad\, a -\mu )^{n-i}(y)
\end{equation}
implies that 
\begin{equation}\label{FFlm}
 F(a, \l , n)\cdot F(a, \mu , m)\subseteq F(a, \l
+\mu , n+m)\;\; {\rm for \; all}\;\; \l , \mu \in K, \;\; n,m\geq
0.
\end{equation}
It follows that the nil-algebra
$$N(a)=F(a,0)=\cup_{n\geq 0}N(a,n), \;\; N(a, n)=F(a, 0,n),$$
is a positively filtered algebra, and that the eigenvalue algebra
$$D(a)=\bigoplus_{\l \in \Ev (\ad \, a)}D(a,\l ), \;\; D(a, \l )=F(a, \l , 0),$$
is  an $\Ev (\ad \, a)$-graded algebra $(D(a, \l )D(a, \mu
)\subseteq D(a, \l +\mu )$ for all $\l , \mu \in \Ev (\ad \, a))$
and the set $\Ev (\ad \, a)$ is an additive submonoid of the field
$K$ since the algebra $A$ is a domain.

For an element $0\neq x\in F(a, \l )$, let $\pi (x)$ be the unique
natural  number $n$ such that $x\in  F(a, \l , n)\backslash F(a,
\l , n-1)$. If $0\neq y \in F(a, \mu )$ then 
\begin{equation}\label{pixy}
\pi (xy)=\pi (x)+\pi (y)
\end{equation}
since by (\ref{FFlm}), $\pi (xy)\leq \pi (x)+\pi (y)$ and,
 by (\ref{dlmn}),
$$(\ad\, a -(\l +\mu ))^{\pi (x)+\pi (y)} (xy)={\pi (x)+\pi (y)\choose \pi
(x)} (\ad\, a - \l)^{\pi (x)} (x)\cdot (\ad\, a -\mu )^{\pi (y)}
(y)\neq 0,$$ since $A$ is a domain and char$\, K=0$.
\begin{lemma}\label{BFD}
Let $A$ be a noncommutative domain, and let $a\in A\backslash
Z(A)$.
\begin{enumerate}
\item $K[a]\subseteq C(a)$, and so $\GK (C(a))\geq  1$.
\item Suppose that $F(a, \l )\neq D(a, \l )$ for some $0\neq \l \in
K$. Then the Gelfand-Kirillov dimension of the  $K$-algebra
$B:=\oplus_{n\geq 0}F(a, n\l )$ is not less than $3$.
\item If $F(a)\neq C(a)$ then $\GK (F(a))\geq 2$.
\end{enumerate}
\end{lemma}

{\it Proof}. The first  statement is evident since the algebra $A$
is a domain and $K=\bK $.

$2$. Choose an element $x\in  F(a, \l , 1)\backslash D(a, \l )$
then $0\neq y:=(\ad \, a-\l )x\in D(a, \l )$.   Note that
$a^iy^jx^k\in F(a, (j+k)\l  , k )\backslash F(a, (j+k)\l , k-1)$
for all $i,j,k\geq 0$  with $j+k>0$. This implies that
$$ B\supseteq \sum_{i,j,k\geq 0}Ka^iy^jx^k=\bigoplus_{i,j,k\geq
0}Ka^iy^jx^k,$$ hence $\GK (B)\geq 3$.

$3$. If $F(a)\neq  C(a)$ then either $N(a)\neq C(a)$ or $D(a) \neq
C(a)$ (or both). In the first (resp.  second) case, we can choose
an element $x\in N(a, 1)\backslash C(a)$ (resp. $0\neq y\in  D(a,
\l )$ for some $0\neq \l \in \Ev (\ad \, a)$). Then
$$ F(a) \supseteq \sum_{i,j\geq 0}Ka^ix^j=\bigoplus_{i,j\geq
0}Ka^ix^j\;\; ({\rm resp. }\;\; F(a) \supseteq \sum_{i,j\geq
0}Ka^iy^j=\bigoplus_{i,j\geq 0}Ka^iy^j).$$ In both cases, this
implies that $\GK (F(a))\geq 2$. $\Box $

\begin{corollary}\label{NorD}
Let a $K$-algebra $A$ be a noncommutative domain of
Gelfand-Kirillov dimension $<3$, and let $a\in A\backslash Z(A)$.
\begin{enumerate}
\item $F(a,\l )=D(a, \l )$ for each $\l  \in K$.
\item Either $F(a)=N(a)$ or $F(a) =D(a)$.
\item The algebra $A$ admits the Dixmier partition.
\end{enumerate}
\end{corollary}

{\it Proof}. $1$. This  follows from Lemma \ref{BFD} and the fact
that $\GK (A)<3$.

$2$. By statement $1$,
$$ F(a)=N(a)\bigoplus\, ( \bigoplus_{0\neq \l \in \Ev(\ad \,
a)}D(a, \l ))=N(a) +D(a).$$ If $F(a)\neq N(a)$ and $F(a) \neq
D(a)$ then $N(a,  1)\neq N(a, 0)$ and $D(a, \l )\neq  0$ for some
$0\neq \l \in \Ev (\ad \, a)$. The algebra $A$ is a domain, so
$0\neq N(a, 1)D(a, \l )$, and so $F(a, \l , 1)\neq F(a, \l , 0)$
by  (\ref{pixy}) which  contradicts to the first statement.

$3$. The statements $2$ and $3$ are equivalent. $\Box $

Suppose that the field $K$ is not necessarily algebraically closed
but still has characteristic zero. Let a $K$-algebra $A$ be a
noncommutative domain. Via the $K$-algebra  monomorphism $A\ra \bA
:=\bK \t A$, $a\ra 1\t a$, we  can identify the algebra $A$ with
its image in $\bA $. The next result  follows from  Corollary
\ref{NorD} (since the type of elements remain the same under field
extension).

\begin{corollary}\label{2NorD}
Let $K$ be a field of characteristic zero (not necessarily
algebraically closed), and let a $K$-algebra $A$ be a
noncommutative domain of Gelfand-Kirillov dimension $<3$ such that
the algebra $\bA =\bK \t A$ is a domain. Then the algebra $A$
admits the Dixmier partition.  $\Box $
\end{corollary}


\section{ Algebras that Satisfy the Commutative Centralizer
Condition, and Dixmier's Problem 6}\label{ccc}

{\sc Definition}.  A noncommutative algebra $A$ satisfies the {\bf
commutative centralizer condition} (the $ccc$, for short) if
 the  centralizer of each element of the
set $A\backslash Z(A)$ is a {\em commutative} algebra.

Clearly, every noncommutative   subalgebra  of an algebra that
satisfies the commutative centralizer condition has this property
as well.

{\em Example 1}. The Weyl algebra $A_1$ satisfies the commutative
centralizer condition \cite{Ami}.

{\em Example 2}. The {\em formal pseudo-differential operator
ring} $L((t^{-1};\d ))$ with coefficients from a field extension
$L/K$ satisfies the commutative centralizer condition where a
derivation $\d $ of $L$ has $\ker \, \d =K$ (Lemma \ref{DSccr}).

{\em Example 3}. The {\em skew Laurent series algebra}
$L=K(H)[[X,X^{-1}; \s ]]$, $\s (H)=H-1$, satisfies the isotropic
centralizer condition. The Weyl algebra $A_1$, its ring of
quotients (so-called, the {\em first Weyl skew field}), and the
noncommutative deformations of type-A Kleinian singularity are
subalgebras of $L$, and so they satisfy the isotropic centralizer
condition (Section \ref{WDRccc}).

{\em Example 4}. The {\em skew Laurent series algebra}
$R=K(H)[[X,X^{-1}; \tau ]]$, $\tau (H)=\l H$, $ \l \in K^*$ is not
a root of $1$, satisfies the commutative centralizer condition.
The quantum plane $\L $ and the quantum  Weyl algebra $A_1(\l )$,
and their rings of quotients  are subalgebras of $R$, and so they
satisfy the commutative centralizer condition (Section
\ref{WDRccc}).

\begin{lemma}\label{PAmaxcomsubr}
 Let $C$ be a commutative subalgebra of an algebra $A$
that satisfies the commutative centralizer condition. The
following statements are equivalent.
\begin{enumerate}
\item $C$ is a maximal commutative subalgebra of $A$.
\item There exists an element $a\in A\backslash Z(A)$ such that $C=C(a)$.
\item $C\neq Z(A)$ and $C=C(b)$, for all elements $b\in C\backslash Z(A)$.
\end{enumerate}
So, for each element $a\in A\backslash Z(A)$, its centralizer
$C(a)$ is the unique maximal commutative subalgebra that contains
the element  $a$.
\end{lemma}

{\it Proof}. $(1\Rightarrow 3)$ Let $C$ be a maximal commutative
subalgebra   of $A$. Then $Z(A)\subseteq C$, but $Z(A)\neq C$
since $A$ is a noncommutative algebra.
 The  centralizer $C(b)$ of each element $b\in C\backslash Z(A)$
 contains $C$. Then $C=C(b)$ since $C(b)$ is an commutative algebra and $C$
 is a maximal commutative subalgebra of $A$.

 $(3\Rightarrow 2)$ Evident.

$(2\Rightarrow 1)$ Let $C=C(a)$, for some $a\in A\backslash Z(A)$,
and let $C'$ be a commutative subalgebra of $A$ containing $C$.
Then $C'\subseteq C$ since $a\in C'$ and $[a, C']=0$. The algebra
$A$ satisfies the $ccc$, hence $C$ is the maximal commutative
subalgebra of $A$. $\Box $

So, Lemma \ref{PAmaxcomsubr} states that for an algebra that
satisfies the {\em ccc} the concepts of centralizer and maximal
commutative subalgebra coincide.


\begin{corollary}\label{PACf=Cg}
 Let an  algebra $A$ satisfy the
commutative centralizer condition.
\begin{enumerate}
\item Let $a,b\in A\backslash Z(A)$. Then $ab=ba$ iff
$C(a)=C(b)$. $ab\neq ba$ iff $C(a)\cap C(b)=Z(A)$.
\item Let $C$ be a maximal commutative subalgebra of $A$, $a\in
A\backslash Z(A)$, and $p(t)\in Z(A)[t]$ be a polynomial of
positive degree with coefficients from $Z(A)$. If $p(a)\in
C\backslash Z(A)$ then $a\in C$.
\end{enumerate}
\end{corollary}

{\it Proof}. $1$. If $C(a)=C(b)$ then $[a,b] =0$ since $A$
satisfies the $ccc$. Suppose that the elements $a$ and $b$
commute. By Lemma \ref{PAmaxcomsubr}.(2),
 $C(a)$ is a maximal commutative subalgebra of $A$. Since $b\in C(a)\backslash
 Z(A)$, by Lemma \ref{PAmaxcomsubr}.(3),
 $C(a)=C(b)$. This proves the first part of the first statement.

If $ab\neq ba$ then $C(a)\cap C(b)=Z(A)$ since otherwise
 we can choose an element, say $c$, from the set $(C(a)\cap C(b))\backslash Z(A)$.
 By Lemma \ref{PAmaxcomsubr}.(3), $C(a)=C(b)=C(c)$, a contradiction.

If $C(a)\cap C(b)=Z(A)$ then $ab\neq ba$ since otherwise by the
first statement $C(a)=C(b)$, hence $C(a)=C(a)\cap C(b)=Z(A)$ which
contradicts to the fact that $a\not\in Z(A)$ $(a\in C(a))$.

$2$. By Lemma \ref{PAmaxcomsubr}.(3), $C=C(b)$ where $b=p(a)$. The
elements $a$ and $b$ satisfy the condition of statement 1, so
$C(a)=C(b)$, hence $a\in C$. $\Box $

\begin{lemma}\label{hKhPZA}
Let an algebra $A$  satisfy the commutative centralizer condition,
and let $a$, $b$, and $p$ be nonzero elements of $A$ such that
$ab=ba$, $[ a,p]=\l p$, and $[ b,p] = \mu p$ for some $\l , \mu
\in K^*$. Then $b\in K^*a+Z(A)$.
\end{lemma}

{\it Proof}. The assumptions $[ a,p]=\l p\neq 0$ and $[ b,p] = \mu
p\neq 0$ imply that $a,b, p\not\in Z(A)$. Since the elements $a$
and $b$ commute and do not belong to the centre of the algebra
$A$, we have $C:=C(a)=C(b)$ and $C\cap C(p)=Z(A)$ (Corollary
\ref{PACf=Cg}). It follows  from $[ \mu a-\l b, p]=0$ that $\mu
a-\l b\in C\cap C(p)=Z(A)$, hence $b\in K^*a+Z(A)$. $\Box $

Let $A$ be a noncommutative algebra. Suppose that the centralizer
$C(a)$ of an
 element $a\in A$ is a commutative algebra. Then $a\not\in Z(A)$ since
$A$ is a noncommutative algebra. Suppose that the set $D(a,\l ,
A)\neq 0$, for some $0\neq \l \in K$, contains an invertible
element, say  $p$, then
  $[a,p]=\l p$. This equality
 is equivalent to the equality $pa=(a-\l )p$. Then $[a,p^{-1}]=-\l p^{-1}$, and  so
$D(a, -\l ,A)\neq 0$. The inner automorphism $\o_p:u\ra pup^{-1}$
of the algebra $A$ preserves the
 centralizer $C(a)$, that is $\o_p (C(a))\subseteq C(a)$, since, for all $c\in C(a)$,
$$ pcp^{-1}a=pc(a+\l )p^{-1}=p(a+\l )cp^{-1}=(a-\l +\l )pcp^{-1}=a\, pcp^{-1}.$$
In the argument  above we use the commutativity  of the
centralizer $C(a)$. The inner automorphism
 $\o_{p^{-1}}$ also preserves the algebra $C(a)$, thus
$$ \o_p (C(a))=C(a).$$
If $q$ is another nonzero element of $D(a,\l ,A)$ then
 $qp^{-1}\in C(a)$ since
$$qp^{-1}a=q(a+\l )p^{-1}=(a-\l +\l )qp^{-1}=aqp^{-1}.$$
This proves that 
\begin{equation}\label{D(a,l)=Cp=pC}
  D(a, \l ,A)=C(a)p=pC(a).
\end{equation}
We denote by $\tau_\l $ the restriction of the inner automorphism
$\o_p$ of the algebra $A$ to the subalgebra $C(a)$. Since the
algebra $C(a)$ is commutative, 
\begin{equation}\label{pc=tcp}
 qc=\tau_\l (c)q, \;\; {\rm for \;   all }\;\; 0\neq q\in D(a, \l
,A), \;\; c\in C(a).
\end{equation}

\begin{proposition}\label{Cafrp}
Suppose that an algebra $A$ satisfies the
 commutative centralizer condition. Let $a$, $b$, $p$, and $q$ be
 nonzero elements of the algebra $A$ such that $ab=ba$,
 $[a,p]=\l p$, and $[b,q]=\mu q$ for some $\l , \mu \in K^*$. If
 the elements $p$ and $q$ are invertible in $A$,  the
 centralizer $C(a)$ of the element $a$ is a domain, and none of
 the elements $p^i$, $i\geq 1$, is an eigenvector for $\ad \, b$
 then the sum
  $M:=\sum_{i,j\geq 0}\, C(a)q^ip^j$ is a free left
$C(a)$-submodule
 of $A$ with the free basis $\{ q^ip^j\, | \, i, j\geq 0\} $.
\end{proposition}

{\it Proof}. The elements $a$ and $b$ commute and do not belong to
the centre of the algebra $A$ since $[a,p]=\l p\neq 0$ and
$[b,q]=\mu q\neq 0$. The algebra $A$ satisfies the $ccc$,  so
$C:=C(a)=C(b)$ is a commutative subalgebra of $A$ which is
  a domain (by the assumption). The elements $p$ and $q$ are units in $A$, so, by
 (\ref{D(a,l)=Cp=pC}),  $D(a,\l ,A)=Cp=pC $ and $D(b,\mu , A)=Cq=qC$. We denote by
$\tau_\l $ and $\tau_\mu $ the restrictions of the inner
automorphisms $\o_p$ and
 $\o_q$ of the algebra $A$ to the subalgebra $C$ as defined in (\ref{pc=tcp}). The element $q$
 is a unit in $A$, so $0\neq q^i\in D(b,i\mu ,A)$, for all $i\geq 0$, and $\{ i\mu \} $ are distinct
 eigenvalues of the inner derivation $\ad\, b$ of the algebra $A$ since the field $K$ has
characteristic zero. Thus

\begin{equation}\label{N=Cqi=qiC}
  N:=\sum_{i\geq 0}\, Cq^i=\bigoplus_{i\geq 0}\,
Cq^i=\bigoplus_{i\geq 0}\, q^iC
\end{equation}
is the free left and right $C$-module since $q^iC=Cq^i$. Moreover,
$q^ic=\tau_\mu^i(c)q^i$, for all $c\in C$.

For an element $u=\sum \, c_iq^i$ of $N$ with $c_i\in C$, the set
$\supp \, u:= \{ i \, | \, c_i\neq 0\} $ is called the {\em
support} of $u$. If $i\in \supp \, u$ then
$$ |\supp (\ad\, b - i\mu )(u)|=|\supp \, u|-1.$$
Suppose that the set $\{ q^ip^j\, |\, i,j\geq 0\} $ is not a free
basis for the left $C$-module $N$.
 Then there exists a nontrivial relation,
$$ r:= u_dp^d+\cdots +u_1p+u_0=0 \;\; {\rm with }\;\; u_i \in N,$$
where $d\geq 1$ and $u_d\neq 0$ (see the definition of the
$C$-module $N$). The number $d$ is called the {\em degree}
 of the relation $r$. We may assume that $d\geq 1$ is the least possible
  degree for all nontrivial relations, then
 $u_0\neq 0$, since otherwise we can divide the relation $r$ by $p$ on the right and
obtain a nontrivial relation of degree $d-1$, a  contradiction. We
may also assume
 that the number $s=|\supp \, u_0|$ is the least possible for all nontrivial
 relations of degree $d$. The element $b$ belongs to the algebra $C$, so, for each $j\geq 1$,
$$ (\ad \, b) (p^j)=bp^j-p^jb=(b-\tau^j_\l (b))p^j=b_jp^j,$$
 where $b_j:=b-\tau^j_\l (b)\in C$. Since $(\ad \, a) (p^j)=j\l
 p^j$, and, by the assumption, for each $j\geq 1$, the element $p^j$
 is not an eigenvector for $\ad \, b$, we conclude that each element
 $b_j$ {\em does not} belong to the field $K$.  Let us fix $k\in \supp \, u_0$, then
 $$ |\supp (\ad \, b -k\mu )(u_0)|<|\supp\, u_0|,$$
 and so the relation
$$ (\ad \, b -k\mu )(r)=\sum^d_{j=1}\, ((\ad\, b -k\mu )(u_j)+u_jb_j)p^j+
 (\ad \, b -k\mu )(u_0)=0$$
must be trivial by the minimality of $d$ and $s$. This means  that
for each $j$ the  coefficient of $p^j$ must be zero. In
particular,
 for  $j=d$, we have
$$ (\ad \, b -k\mu )(u_d)+u_db_d=0.$$
 The element $u_d$ has the form $\sum \, c_iq^i$ where not all $c_i\in C$ are equal to zero.
Substituting the sum $u_d=\sum\, c_iq^i$ in the equality above and
using
 $(\ad\, b)(q^i)=i\mu q^i$, we obtain the equality
$$ 0=\sum ((i-k)\mu c_iq^i+c_iq^ib_d)=\sum c_i((i-k)\mu +\tau^i_\mu (b_d))q^i.$$
 By (\ref{N=Cqi=qiC}), for each $i$, the coefficient of $q^i$ must be $0$. Since not all elements $c_i$ are zero  and
 $C$ is a domain (by the assumption), we must have, for some $i$ (such that $c_i\neq 0)$,
$$ (i-k)\mu +\tau^i_\mu (b_d)=0,$$
 hence $b_d\in K$, a contradiction. Thus $M$ is a free left $C$-module with  the
basis $\{ q^ip^j\} $. $\Box $

Recall that commuting non-central elements of an algebra that
satisfies the {\em ccc} have common centralizer (Corollary
\ref{PACf=Cg}.(1)). The next result shows that they also have
common nil-algebra.

\begin{proposition}\label{NAappb}
Suppose that an algebra $A$ satisfies the commutative centralizer
condition, and that elements $a,b\in A\backslash Z(A)$ commute.
Then $N(a,A)=N(b,A)$, and  $N(a,n,A)=N(b,n,A)$ for all $n\geq 0$.
\end{proposition}

{\it Proof}. It suffices to prove that $N(a,n,A)=N(b,n,A)$ for all
$n\geq 0$. We use induction on $n$. The algebra $A$ satisfies the
commutative centralizer condition, and the non-central elements
$a$ and $b$ commute, hence $C(a)=C(b)$ (Corollary \ref{PACf=Cg}),
and the case $n=0$ holds.
 The derivations $\der :=\ad \, a$ and $\d :=\ad \, b $ commute
 ($\der \d -\d \der =\ad \, [ a,b]=\ad \, 0=0$).
 Suppose that $N(a,i,A)=N(b,i,A)$ for
 all $i\leq n$. Then
\begin{eqnarray*}
 f\in N(a, n+1,A)&\Leftrightarrow& 0=\der^{n+2}f=\der^{n+1}\der f
 \Leftrightarrow \der f\in N(a, n,A)=N(b,n,A)\Leftrightarrow\\
0=\d^{n+1}\der f = \der \d^{n+1}f&\Leftrightarrow& \d^{n+1}f\in
C(a,A)=C(b,A)\Leftrightarrow
0=\d \d^{n+1}f=\d^{n+2}f,\\
\end{eqnarray*}
thus $N(a, n+1,A)=N(b,n+1,A)$. By induction, $N(a,n,A)=N(b,n,A)$
for all $n\geq 0$.  $\Box $

{\bf Proof of Theorem \ref{G3Apr}}.

The $K$-algebra $A$ has Gelfand-Kirillov dimension $<3$, and $\bA
$ is a domain, so $A$ admits the Dixmier partition (Corollary
\ref{2NorD}). The algebra $A$ satisfies the {\em ccc} as a
noncommutative subalgebra of the algebra $\G $ that satisfies the
{\em ccc}. So, it remains to show that the algebra $A$ satisfies
the {\em hcc}. Let $a\in A\backslash Z(A)$, and let $C=C(a)$.

 Suppose that the algebra $C$ contains a
strongly nilpotent element (resp. a weakly nilpotent element).
Then, by Proposition \ref{NAappb}, all the elements of
$C\backslash Z(A)$ are strongly nilpotent (resp. weakly
nilpotent).

Suppose that $C$ doest not contain a nilpotent element. Then
either all the elements of $C\backslash Z(A)$ are generic or,
otherwise, there exists a semi-simple element, say $a\in
C\backslash Z(A)$. Note that the field extension does not change
the type of elements, $\GK_{\bK }(\bA )=\GK (A)<3$,  $Z(\bA )=\bK
\t Z(A)$, and $\bA $ is a domain. So, first let us  assume that
$K=\bK $. Suppose that $b$ is another semi-simple element of $C$.
Then $ab=ba$ ($C$ is a commutative algebra  since the algebra  $A$
satisfies the {\em ccc}), $[a,p]=\l p$, and $[b,q]=\mu q$ for some
non-zero elements $p,q\in A$ and some nonzero  scalars $\l ,\mu
\in K^*$. We claim that, for some $i\geq 1$, the element $p^i$ is
an eigenvector for  the inner derivation $\ad \, b$  of the
algebra $A$. Suppose that  this is not true, then, by  Proposition
\ref{Cafrp}, the sum $M=\sum_{i,j\geq 0}Cq^ip^j$ is a free  left
$C$-submodule of $A$  with the free basis $\{ q^ip^j\, | \,
i,j\geq 0\}$. The element $a$ is not an algebraic element (since
$a\not\in K$ and $K=\bK$), so the subalgebra $K[a]$ of $C$
generated by the element $a$ is isomorphic to a polynomial algebra
in one indeterminate. Let $\L $ be the $K$-subalgebra
 of  $A$  generated by the elements $a$, $p$, and $q$, and
let $\{ \L_n\}_{n\geq 0}$ be the standard filtration of the
algebra $\L $ determined by the total degree of the
(non-commutative) generators. For each $n\geq 0$,
 $$ \L_n \supseteq \bigoplus_{i,j,k\geq 0}\, \{ Ka^iq^jp^k\, | \, i+j+k\leq n\}, $$
 and so
 $$ \dim_K \,\L_n\geq {n+3\choose 3}={(n+1)(n+2)(n+3)\over 3!},
 $$
 hence the Gelfand-Kirillov dimension $\GK (\L )\geq 3$. On the
  other hand, $\L $ is a subalgebra of the algebra $A$,
   and so $ \GK (\L ) \leq \GK (A)<3$, a  contradiction. We have
  proved that, for some $i\geq 1$, the element $p^i$ is a common
  eigenvector for the inner derivations $\ad \, a$ and $\ad \, b$:
  $[a,p^i]=i\l p^i\neq 0 $ and $[b, p^i]=\g p^i$ for some $\g \in
  K$. We must have $\g \neq 0$ since otherwise the non-central
  elements $b$ and $p^i$ commute, and so they have the common
  centralizer $C(b)=C(p^i)$. By the assumption, the non-central
  elements $a$ and $b$ commute, so $C(a)=C(b)$ which implies that
  the elements $a$ and $p^i$ commute which is impossible.
 Now, by Lemma \ref{hKhPZA}, $b\in K^*a+Z(A )$, hence $K^*a+Z(A)$
 is the set of all semi-simple elements of $C$ (where $K=\bK $).

 Suppose that $K\neq \bK $. Since $a\in A$, $C(a, \bA )=\bK \t C(a,
 A)$, $Z(\bA )=\bK \t Z(A)$, and $\bK^*a+Z(\bA )$ is the set of all the
 semi-simple elements in $C(a,\bA )$, we see that
 $K^*a+Z(A )$ is the set of all the
semi-simple elements in $ C(a, A)$. So, the algebra $A$ satisfies
the {\em hcc}.  This finishes the proof of the theorem. $\Box $

\begin{theorem}\label{DP6ccc}
{\sc (An Analog of Dixmier's Problem 6 for Algebras that satisfy
the Commutative Centralizer Condition)} Let $B=\cup_{i\geq 0}B_i$
be a filtered algebra over an algebraically closed field of
characteristic zero such that the associated graded algebra $\gr
\, B$ is a commutative domain, and $B$ satisfies the commutative
centralizer condition and admits the Dixmier partition. Let
$f(t)\in K[t] $ be an arbitrary polynomial  of degree $>1$.
\begin{enumerate}
\item Let  $i=1,2 $. If $a\in \D_i(B)$ then $f(a)\in \D_i(B)$ provided $f(a)\not\in Z(B)$.
\item Let $i=3,4$. If $a\in \D_i(B)$  then $f(a)\in \D_5(B)$.
\item  If $a\in \D_5(B)$ then $f(a)\in \D_5(B)$.
\end{enumerate}
\end{theorem}

{\it Proof}. By Corollary \ref{cor11}, $a\not\in Z(B)$ and
$f(a)\in Z(B)$ imply $a\in B_0$. Clearly $B_0\subseteq Z(B)\cup
\D_{1,2}(B)$ since $\gr (B)$ is a commutative algebra.

 Let $i=1,2$. By Proposition \ref{NAappb}, $a\in
\D_i(B)\Leftrightarrow  f(a)\in \D_i(B)$  provided $f(a)\not\in
Z(B)$.

For $i=3,4$, $a\not\in B_0$, hence $f(a)\not\in Z(B)$,
$f(a)\not\in \D_{3,4}(B)$, by Theorem \ref{Beigvfa=0}, and
$f(a)\not\in \D_{1,2}(B)$, by the argument above. Therefore,
$f(a)\in \D_5(B)$.

For $i=5$, $a\not\in B_0$, hence $f(a)\not\in Z(B)$, $f(a)\not\in
\D_{3,4}(B)$ (Theorem \ref{Beigvfa=0}),
 and $f(a)\not\in \D_{1,2}(B)$, hence $f(a)\in \D_5(B)$. $\Box $


\section{The Weyl Division Ring satisfies the Commutative
Centralizer Condition}\label{WDRccc}

Let $L$ be a field, and let $\s $ be an automorphism of the field
$L$. The {\em skew Laurent series algebra} $\CL =L[[X, X^{-1}; \s
]]$ consists of (inverse) Laurent series
$a=\sum_{i=-\infty}^na_iX^i$ with $n\in \mathbb{Z}$ and
 coefficients  $a_i\in L$. Addition and multiplication in $\CL $ are given by
 the formulae:
$$ a+b= \sum \, (a_i+b_i)X^i \;\; {\rm and } \;\;
ab=\sum_i(\sum_{j+k=i}\, a_j\s^j (b_k))X^i,$$ where $b=\sum \,
b_iX^i$. The ring $\CL $ is a division ring. If $a\neq 0$ then the
{\em degree} $\deg \,a$ of $a$  is the maximal integer $n$ such
that $a_n\neq 0$, and so $a=a_nt^n+\cdots $ where $a_nt^n$ (resp.
$a_n$) is the {\em leading term} (resp. the {\em leading
coefficient}) of $a$, and three dots denote terms of smaller
degree. We define $\deg (0)=-\infty $. For $a=a_nX^n+\cdots
,b=b_mX^m+\cdots \in \CL $, $ab=a_n\s^n(b_m)X^{n+m}+\cdots $,  and
so $\deg (ab)=\deg (a)+\deg (b)$ and $\deg (a+b)\leq \max \{ \deg
(a),\deg (b)\}$.

The subring $L[X,X^{-1};\s ]$ of $\CL$ that consists of all finite
sums $\sum\, a_iX^i$ is called the {\em skew Laurent polynomial
ring}. The ring $\CL $ is a completion of $L[X,X^{-1};\s ]$ with
respect to the valuation $v(a):=-\deg \, a$.

The subring $\CL_-:=\{ a\in \CL \, |\, \deg \, a\leq 0 \}$ of $\CL
$ is a {\em skew series algebra}
 $L[[ X^{-1}; \s^{-1}]]$. It is a local ring with  (unique) maximal ideal
 $\CL_-X^{-1}$. The residue ring $\CL_-/\CL_-X^{-1}$ is canonically isomorphic to
 the field $L$. A nonzero element of $\CL_-$ is a unit iff it has degree $0$.

 For a field extension $E$ of $F$ we denote by $[E:F]$ the
 dimension of $E$ over $F$. Let $L^\s :=\{ l\in L\, | \, \s
 (l)=l\}$ be the subfield of $L$ of fixed elements for the
 automorphism $\s $ of the field $L$.

\begin{theorem}\label{CLccc}
The skew Laurent series ring $\CL =L[[X, X^{-1}; \s ]]$ with
coefficients from a field $L$ satisfies the commutative
centralizer condition provided $\s \neq id_L$ and, for each $l\in
L\backslash L^\s $, $\si (l)\neq l$ for all $i\geq 1$.
\begin{enumerate}
\item If $u\in L\backslash L^\s$, then $C(u)=L$.
\item If $u\in \CL \backslash L$, $\deg \, u=0$, and the leading coefficient
$u_0$ of $u$ does not belong to $L^\s $, then $C(u)\subseteq
\CL_-$, and the map $C(u)\ra L=\CL_-/\CL_-X^{-1}$, $c\ra
c+\CL_-X^{-1}$, is a $L^\s$-algebra isomorphism, $[C(u):L^\s
(u)]=[L:L^\s (u_0)]$.
\item If $u\in \CL \backslash L$ with $\deg\, u\neq 0$ then
 $C(u)=L^\s [[v, v^{-1}]]$ is a Laurent series field with coefficients from
 the field $L^\s $  where $v$ is an element of $C(u)$ which has
 the least positive degree.
 $[\, C(u):L^\s [[u, u^{-1}]] \, ]=|{\deg \, u\over \deg \, v}|<\infty $.
\end{enumerate}
\end{theorem}

{\it Proof}. Observe that the centralizer $C=C(u)$ of any nonzero
element of $\CL $ is a skew subfield of $\CL $.

$1$. By the assumption $u\not\in  L^\s =L^{\si}$ (for all $0\neq
i\in \mathbb{Z}$), which means that $\si (u)-u\neq 0$ for all
$0\neq i\in \mathbb{Z}$. Now, $a=\sum \, a_i X^i\in C(u)$ iff
 $0=au-ua= \sum \, a_i(\si (u)-u)X^i$ iff $a_i=0$ for all $0\neq
i\in \mathbb{Z}$ iff $C(u)=L$.

$2$. By the assumption, the element $u$ has the form $\sum_{j\leq
0}\, u_jX^j$ with  $u_0\in L\backslash L^\s$. Let
 $0\neq a=\sum \, a_iX^i \in C(u)$. Then the leading term of $a$, say $a_dX^d$,
 commutes with $u_0$, hence $d=0$ since $u_0\not\in L^{\si}=L^\s $ for all
 $0\neq i\in \mathbb{Z}$. Thus  $a=\sum_{i\leq 0} \, a_iX^i $ and
$C\subseteq \CL_-$. Now,
$$ 0=au-ua=\sum_{s\leq 0}\, (\sum_{i+j=s}\, (a_i\si (u_j)-u_j\s^j (a_i)))X^s.$$
Equating the coefficients of $X^s$ to $0$ in the equality above we
obtain the infinite
 system of equations with unknowns $a_i$, $i\leq 0$,
$$
 (\s^s (u_0)-u_0)a_s+\sum_{i+j=s,\,  i>s}\, (a_i\si (u_j)-u_j\s^j (a_i))=0,
\; s\leq 0.$$ For a given $a_0$, the system has a unique solution
since $\s^s (u_0)-u_0\neq 0$,
 for all $s< 0$ (since $u_0\not\in L^\s$). So, the map
$C(u)\ra K(H)$, $c\ra c+\CL_-X^{-1}$, is a $L^\s$-algebra
isomorphism, hence
 $[C(u):L^\s (u)]=[L:L^\s (u_0)]$.

$3$. Since $C(u)=C(u^{-1})$ we may assume that $d:=\deg \, u>0$.
Then the set $ G:=\{ \deg \, c \, | \,  \; 0\neq c\in C(u)\}$ is a
nonzero subgroup of $\mathbb{Z}$ that contains $ \mathbb{Z}d$,
hence $G=\mathbb{Z}t$ where an integer $t$ is the least positive
 element of $G$. Clearly, $t|d$. Fix an element $v\in C(u)$ with $\deg \,
 v=t$. We claim that $C(u)=L^\s [[v,v^{-1}]]$. Obviously,
 $L^\s [[v,v^{-1}]]\subseteq C(u)$. In order to prove the reverse
 inclusion, we first show that, for a given $\deg \, c$, the
 leading term of $c\in C(u)$ is unique up to  a factor of $L^\s$.

 So, let elements $c,c'\in C(u)$ have the same degree, say $n$,
 and the leading terms $\alpha X^n$ and  $\alpha' X^n$
 respectively. These leading terms commute with the leading term,
 say $u_dX^d$, of the element $u$:
 $$ 0=[\alpha X^n, u_dX^d]=(\alpha \s^n(u_d)-u_d\s^d(\alpha
 ))X^{n+d} \; {\rm and}\; 0=[\alpha' X^n, u_dX^d]=(\alpha'
 \s^n(u_d)-u_d\s^d(\alpha'
 ))X^{n+d},$$
and so $ \frac{\s^d(\alpha )}{\alpha }=\frac{\s^n(u_d )}{u_d
}=\frac{\s^d(\alpha' )}{\alpha' }$ which implies that
$\s^d(\frac{\alpha}{\alpha'})=\frac{\alpha}{\alpha'}$, hence
$\frac{\alpha}{\alpha'}\in L^{\s^d}=L^\s$, as required.

 Let $a$ be an arbitrary element of $C(u)$. The elements $a$ and
 $v^{i_1}$, $i_1:=t^{-1}\deg \, a$, of $C(u)$ have the same degree
 $\deg \, a$, so we can choose a scalar, say $\l_1\in L^\s $, such
 that the elements $a$ and $\l_1v^{i_1}$ have the same leading
 term. Hence the element
 $a_1:=a-\l_1v^{i_1}\in C(u)$ has degree $\deg \, a_1<\deg \, a$. Repeating the same argument
 for the element $a_1\in C(u)$, we can find a scalar $\l_2\in L^\s$
 and an integer
 $i_2 \in \mathbb{Z}$ such that the element $a_2:=a-\l_2v^{i_2}\in C(u)$
 has  degree $\deg \, a_2<\deg \, a_1$. Proceeding in this way, we can find infinitely
many elements of $C(u)$,
$$a_\nu :=a_{\nu -1}-\l_\nu v^{i_\nu } , \; \l_\nu \in L^\s , \; \nu \geq 1,$$
such that $\deg \, a_1>\deg \, a_2>\cdots $ and $i_1>i_2 >\cdots
$.
 So,
$$ a=\sum_{\nu \geq 1}\, \l_\nu v^{i_\nu } \in L^\s [[V,V^{-1}]],$$
thus $C(u)\subseteq L^\s [[V,V^{-1}]]$, as required.  A degree
argument shows that $[C(u,):L^\s [[u, u^{-1}]] ]=|{\deg \, u\over
\deg \, v}|<\infty $. $\Box $

\begin{corollary}\label{cenCL}
Let $\CL $ be as in Theorem \ref{CLccc}. Then the  centre $Z(\CL
)=L^\s $.
\end{corollary}

{\it Proof}. By Theorem \ref{CLccc}.(1,3), $C(X)=L^\s [[X,
X^{-1}]]$ and $C(u)=L$ for each $u\in L\backslash L^\s $. Now,
$$ L^\s \subseteq Z(\CL )\subseteq C(u)\cap C(X)=L\cap L^\s [[X,
X^{-1}]]=L^\s , $$ and so $Z(\CL )=L^\s $. $\Box $

Let $K(H)$ and $K(H,C)$ be fields of rational functions with
coefficients from the field $K$ of characteristic zero in one and
two variables. It can be easily verified that the following four
division rings satisfy the conditions of Theorem \ref{CLccc}, and
so {\em satisfy the commutative centralizer condition}.

\begin{enumerate}
\item $\CA :=K(H)[[X,X^{-1}; \s ]]$, $\s\in \Aut_K(K(H))$, $\s
(H)=H-1$.
\item $\CB_\l :=K(H)[[X,X^{-1};\tau ]]$, $\tau \in \Aut_K(K(H))$,
$\tau (H)=\l H$ where $\l \in K^*$ is not an $i^{th}$ root of $1$
for all $i\geq 1$.
\item $\CC :=K(H,C)[[X,X^{-1}; \s ]]$, $\s\in \Aut_K(K(H,C))$, $\s
(H)=H-1$, $\s (C)=C$.
\item $\CE_\l :=K(H,C)[[X,X^{-1};\tau ]]$, $\tau \in \Aut_K(K(H,C))$,
$\tau (H)=\l H$, $\tau (C)=C$,  where $\l \in K^*$ is not an
$i^{th}$ root of $1$ for all $i\geq 1$.
\end{enumerate}
By Corollary  \ref{cenCL},
$$ Z(\CA )=K, \; Z(\CB_\l )=K, \; Z(\CC )=K(C), \; Z(\CE_\l
)=K(C).$$

Next we show that certain algebras (considered in Introduction)
are subalgebras of one of these division rings, and then applying
Theorem \ref{G3Apr} we prove that they admit the Dixmier partition
and satisfy both the homogeneous centralizer condition and the
commutative centralizer condition.

The Weyl algebra $A_1$ can be identified with its image in the
division ring $\CA $ via the $K$-algebra monomorphism $A_1\ra
\CA$, $x\ra X$, $\der \ra HX^{-1}$, as follows from the commutator
 $(HX^{-1})X-X(HX^{-1})=H-(H-1)XX^{-1}=1$. The Weyl algebra $A_1$
 is a Noetherian domain, its division ring $Q(A_1)$ (so-called,
 the {\em Weyl division ring}) is a subring of $\CA $ by the
 universality of localization, and so the Weyl division ring
 $Q(A_1)$ satisfies the {\em ccc}.

Clearly, the Weyl division ring $Q(A_1)$ has the trivial centre
\begin{equation}\label{ZCL=ZG=K}
Z(Q(A_1) )=K
\end{equation}
since (by Theorem \ref{CLccc}) $C(X, \CA )=K[[X,X^{-1}]]$,
$C(H=\der x, \CA )=K(H)$, and so
\begin{eqnarray*}
K & \subseteq & Z( Q(A_1)) \subseteq C(H, Q(A_1))\cap C(X,Q(A_1))
=Q(A_1) \cap C(H, \CA )\cap C(X, \CA ) \\
&=& Q(A_1)\cap  K(H)\cap K[[X, X^{-1}]]=  Q(A_1) \cap K=K.
\end{eqnarray*}

Given a nonzero polynomial $a=a(H)\in K[H]$. The algebra
 $$ A(a):=K\langle X,Y,H\, | \, XH=(H-1)X, \, YH=(H+1)Y, \, YX=a(H),
 \, XY=a(H-1)\rangle $$
is called the {\em noncommutative deformation of type-A Kleinian
singularities} \cite{bavaia}, \cite{hodges}, \cite{CBH}.
 The algebra $A(a)$ is a subalgebra of the division ring $\CA $
 via the $K$-algebra monomorphism $A(a)\ra \CA $, $X\ra X$, $Y\ra
 aX^{-1}$, $H\ra H$.  Note that all prime infinite dimensional
 factor algebra $B_\l :=Usl(2)/(C-\l )\simeq A(\l -H(H+1))$ are of this
 type where $\l \in K$ and $C$ is the {\em Casimir} element of $Usl(2)$
 (see  Section \ref{Usl2}).

 The Gelfand-Kirillov dimension $\GK (A_1)=\GK (A(a))=2$. So, the
 next result follows from Theorem \ref{G3Apr}, and
 is a generalization
 to the Weyl division ring of the result of Amitsur that the centralizer of an arbitrary
 non-scalar element is a commutative algebra \cite{Ami}.

\begin{corollary}\label{QA1ch}
\begin{enumerate}
\item The Weyl division ring $Q(A_1))$ satisfies the commutative
centralizer condition.
\item The first Weyl algebra $A_1$, the noncommutative deformations of the  type-A
 Kleinian singularities, and any noncommutative subalgebra of the
 Weyl division ring $Q(A_1)$ with Gelfand-Kirillov dimension  $<3$
 admits the Dixmier partition and satisfy both the homogeneous
 centralizer condition and the commutative centralizer condition.
 $\Box $
 \end{enumerate}
\end{corollary}

The {\em quantum plane} $\L =K\langle x,y\,| \, xy=\l yx\rangle $,
where  $\l \in K^*$ is not an $i^{th}$ root of $1$ for all $i\geq
1$,  can be identified with its image in the division ring $\CB_\l
$ under the $K$-algebra monomorphism $\L \ra \CB_\l$, $x\ra X$, $y
\ra HX^{-1}$, as follows from
 $X(HX^{-1})=\l (HX^{-1})X$. The quantum plane $\L $
 is a Noetherian domain of Gelfand-Kirillov dimension $2$, its division ring $Q(\L )$  is a subring of $\CB_\l $ by the
 universality of localization, and so the division ring
 $Q(\L)$ satisfies the {\em ccc}.

The {\em quantum Weyl algebra} $A_1(\mu ):=K\langle x,\der\,| \,
\der x-\mu x\der =1 \rangle $, $1\neq \mu \in K^*$,  is a
subalgebra of the division ring $\CB_\mu $ via the $K$-algebra
monomorphism $A_1(\mu )\ra \CB_\mu $, $\der \ra X$, $x\ra (\mu
-1)^{-1}(Y-1)X^{-1}$,  as follows from
$$ X (\mu
-1)^{-1}(Y-1)X^{-1}-\mu (\mu -1)^{-1}(Y-1)X^{-1} X=(\mu
-1)^{-1}(\mu Y-1 -\mu (Y-1))=1.$$

  The quantum
Weyl algebra $A_1(\mu )$ is a Noetherian domain of
Gelfand-Kirillov dimension $2$, its division ring $Q(A_1(\mu ))$
is a subalgebra of $\CB_\mu $, and so satisfies the {\em ccc}.

\begin{corollary}\label{QLch}
\begin{enumerate}
\item The division rings of the quantum plane and the quantum Weyl algebra
 satisfy the commutative centralizer condition.
\item The quantum plane $\L $, the quantum Weyl algebra $A_1(\mu )$,
 and any noncommutative subalgebra of the
 division ring $Q(\L )$ with Gelfand-Kirillov dimension  $<3$
 admits the Dixmier partition and satisfy both the homogeneous
 centralizer condition and the commutative centralizer condition.
 $\Box $
 \end{enumerate}
\end{corollary}


\section{The Algebras $Usl(2)$ and $U_qsl(2)$ Admit the Dixmier Partition and Satisfy
 the $ccc$ and $hcc$} \label{Usl2}

In this section we prove that the universal enveloping algebra
$Usl(2)$ of the Lie algebra $sl(2)$ and its quantum analog
$U_qsl(2)$ admit the Dixmier partition and satisfy
 the $ccc$ and $hcc$ (Propositions \ref{Usl2p3} and \ref{Uqsl2p})). We can not  apply directly Theorem \ref{G3Apr}
 since these algebras have Gelfand-Kirillov dimension $3$, but
 modifying arguments slightly the proof of the results for the
 algebras $Usl(2)$ and $U_qsl(2)$ proceeds along the line of the
 proof of Theorem \ref{G3Apr}.

The universal enveloping algebra
$$Usl(2)=K\langle X,Y,H\, | \,
[H,X]=X, \; [H,Y]=-Y, \;  [X,Y]=2H\rangle $$ of the Lie algebra
$sl(2)$ is a $K$-subalgebra of the division ring $\CC $ via the
$K$-algebra monomorphism
$$Usl(2)\ra \CC , \; X\ra X, \; Y\ra (C-H(H+1))X^{-1}, \; H\ra H,$$
where $C:=YX+H(H+1) $  is the {\em Casimir element} of the algebra
$Usl(2)$, an algebra generator for the centre $Z(Usl(2))=K[C]$ of
the algebra $Usl(2)$.

The universal enveloping algebra $Usl(2)$ is a Noetherian domain
of Gelfand-Kirillov dimension $3$, so its division ring
$Q(Usl(2))$ is a subalgebra of the division ring $\CC $, and so
both algebras $Usl(2)$ and $Q(Usl(2))$ satisfy the $ccc$. The
centre $Z(Q(Usl(2)))=K(C)$, the field of rational functions, as
follows from (using Theorem \ref{CLccc}) $$K(C)\subseteq
Z(Q(Usl(2)))\subseteq \CC \cap C(H, \CC )\cap C(X, \CC )=\CC \cap
K(H,C)\cap K(C)[[X, X^{-1}]]=K(C).$$ So, we have proved the first
statement of the next result.

\begin{proposition}\label{Usl2p3}
\begin{enumerate}
\item The universal enveloping algebra   $Usl(2)$  and its division ring
$Q(Usl(2))$ satisfy the commutative centralizer condition, and
the centre $Z(Q(Usl(2)))=K(C)$.
\item The universal enveloping algebra   $Usl(2)$, or any
noncommutative subalgebra $A$ of the  division  ring $Q(Usl(2))$
with Gelfand-Kirillov dimension $\GK (A)<4$ and such that $A\cap
K(C)\neq K$ admits the Dixmier partition and satisfies  both the
homogeneous  centralizer condition and the commutative centralizer
condition.
 \end{enumerate}
\end{proposition}

{\it Proof}.  $2$. We split the proof into several steps.

{\it Step 1}. $\GK (C(a,A))\geq 2$ {\em for all } $a\in A
\backslash Z(A)$.

In order to prove this fact, let us fix an element, say $z$, of
the set $(A\cap  K(C))\backslash K$. Then the polynomial algebra
$K[z]$ is the subalgebra of the centre $Z(A)$. Let $a\in
A\backslash Z(A)$. The centralizer $C(a,A)$ is a commutative
algebra (since the division ring $Q(Usl(2))$ satisfies the $ccc$)
that contains the elements $a$ and $z$. It suffices to show that
the elements $a$ and $z$ are algebraically independent.

If $\deg \, a\neq 0$ (in $\CC $) then this fact is evident since
$\deg \, z=0$.

If $a\in K(H,C)$ then $a\not\in K(C)$ (otherwise $a\in  A\cap
K(C)=A\cap Z(Usl(2))\subseteq Z(A)$, a contradiction), and so $a$
is not an algebraic elements over  the field $K(C)$, hence $a$ and
$z$ are algebraically independent.

In the remaining case when $\deg \, a=0$,  $a\not\in  K(H,C)$,
suppose that  the elements $a$ and $z$ are algebraically
dependent. This implies that $a$ is an algebraic element over the
field $K(C)$ since $z\in K(C)$. Then the leading term, say $a_0\in
K(H,C)$, of $a$ is algebraic over the field $K(C)$, and so $a_0\in
K(C)$. Then the element $a-a_0$ is algebraic over $K(C)$ which is
impossible since $\deg (a-a_0)<0$. This contradiction proves Step
$1$.

{\it Step 2.  The algebra $A$ admits the Dixmier  partition}.

The algebra $\overline{\CC }=\bK (H,C)[[X,X^{-1}; \s ]]$ is a
domain, so is $\bA $ since $\bA \subseteq \overline{\CC}$. So,
without loss of generality we may assume that $K=\bK $. Suppose
that the algebra $A$ does not admit the Dixmier partition. Then
there exists an element $a\in A \backslash Z(A)$ such that
$N(a)\neq C(a)$ and $D(a) \neq C(a)$. Let us fix  nonzero elements
$x\in N(a, 1, A)\backslash C(a)$ and $y\in  D(a, \l , A)$ for some
$0\neq \l \in K$. Then, for all $i,j\geq  0$,
$$  C(a)^*x^iy^j\subseteq F(a, j\l , i,A)\backslash F(a, j\l , i-1,
A)$$ where $C(a)^*:=C(a) \backslash \{ 0\} $, and so  the  torsion
algebra
$$ F(a) \supseteq \sum_{i,j\geq 0}C(a)x^iy^j=\bigoplus_{i,j\geq
0}C(a)x^iy^j\supseteq \bigoplus_{i,j\geq 0}K[a,z]x^iy^j,$$ hence
$\GK (F(a))\geq 4$ which  contradicts to the fact that $\GK
(F(a))\leq \GK  (A)<4$. This proves Step $2$.

{\it Step 3. The algebra $A$ satisfies the} $hcc$.

We still may assume that $K=\bK $. Let $a\in A\backslash  Z(A)$,
and let $C=C(a)$. Suppose that the algebra $C$ contains a strongly
nilpotent element (resp. a weakly nilpotent element). Then, by
Proposition \ref{NAappb}, all the elements of the set $C\backslash
Z(A)$ are strongly nilpotent (resp. weakly nilpotent).

Suppose that the algebra $C$ does not contain a nilpotent element.
Then either all the elements of the set $C\backslash Z(A)$ are
generic or, otherwise, there exists a semi-simple element, say
$a\in C\backslash Z(A)$. Suppose that $b$ is another semi-simple
element of $C$. Then $ab=ba$ since $C$ is a commutative algebra,
$[a,p]=\l p$, and $[b,q]=\mu q$ for some nonzero elements $p,q\in
A$ and some nonzero scalars $\l , \mu \in K^*$. Then, for some
$i\geq 1$, the element $p^i$ is an  eigenvector for the inner
derivation $\ad \, b$ of the algebra $A$ since otherwise, by
Proposition \ref{Cafrp}, the sum $M=\sum_{i,j\geq 0}Cq^ip^j$ is a
free left $C$-submodule of $A$ with the free basis $\{ q^ip^j\, |
\, i,j\geq 0\} $. By Step 1, $C$ contains the polynomial   algebra
$K[a,z]$, and so $A\supseteq \oplus_{i,j, k,l\geq
0}Ka^kz^lq^ip^j$. Therefore the subalgebra $\L $ of $A$ generated
by the elements $a$, $z$, $p$, and $q$ has Gelfand-Kirillov
dimension $\geq 4$, which contradicts to the fact that $\GK (\L
)\leq \GK (A)<4$. So, we have proved that, for some $i\geq 1$,
$p^i$ is a common eigenvector for  the inner derivations $\ad \,
a$ and $\ad \, b$: $[a,p^i]= i\l p^i$ and $[b,p^i]=\g p^i$ for
some $\g \in K$. We must have $\g \neq 0$ since  otherwise  the
non-central elements $b$ and $p^i$ commute, and so they have the
common centralizer $C(b)=C(p^i)$. By the assumption, the
non-central elements $a$ and $b$ commute, so $C(a)=C(b)$ which
implies that  the elements $a$ and $p^i$ commute which is
impossible. Now, by Lemma \ref{hKhPZA}, $b\in K^*a+Z(A)$, hence
$K^*a+Z(A)$ is the set of all semi-simple elements of $C$. This
finishes the proof of Step 3 and the proposition.  $\Box $

Suppose that a scalar  $q\in K^*$ is {\em not} a root of $1$, and
let $h:=q-q^{-1}$.  The quantum
$$U_qsl(2)=K\langle X,Y,H, H^{-1}\, | \,
XH=qHX, \; YH=q^{-1}HY, \;  [X,Y]=h^{-1}(H^2-H^{-2})\rangle $$ is
a $K$-subalgebra of the division ring $\CE_q $ via the $K$-algebra
monomorphism
$$U_qsl(2)\ra \CE_q , \; X\ra X, \; Y\ra
(C+(2h)^{-1}(\frac{H^2}{q^2-1}-\frac{H^{-2}}{q^{-2}-1}))X^{-1}, \;
H\ra H,$$ where
$C:=YX-(2h)^{-1}(\frac{H^2}{q^2-1}-\frac{H^{-2}}{q^{-2}-1})  $  is
the (quantum) {\em Casimir element} of the algebra $U_qsl(2)$, an
algebra generator for the centre $Z(U_qsl(2))=K[C]$ of the algebra
$U_qsl(2)$.

The  algebra $U_qsl(2)$ is a Noetherian domain of Gelfand-Kirillov
dimension $3$, so its division ring $Q(U_qsl(2))$ is a subalgebra
of the division ring $\CE_q $ which implies that  both algebra
$U_qsl(2)$ and $Q(U_qsl(2))$ satisfy the $ccc$. The centre
$Z(Q(U_qsl(2)))=K(C)$ since $$K(C)\subseteq
Z(Q(U_qsl(2)))\subseteq \CE_q \cap C(H, \CE_q )\cap C(X, \CE_q
)=\CE_q \cap K(H,C)\cap K(C)[[X, X^{-1}]]=K(C).$$ So, we have
proved the first statement of the next result.

\begin{proposition}\label{Uqsl2p}
\begin{enumerate}
\item The quantum $U_qsl(2)$  and its division ring
$Q(Us_ql(2))$ satisfy the commutative centralizer condition, and
the centre $Z(Q(U_qsl(2)))=K(C)$.
\item The quantum   $U_qsl(2)$, or any
noncommutative subalgebra $A$ of the  division  ring $Q(U_qsl(2))$
with Gelfand-Kirillov dimension $\GK (A)<4$ and such that $A\cap
K(C)\neq K$ admits the Dixmier partition and satisfies  both the
homogeneous  centralizer condition and the commutative centralizer
condition.
 \end{enumerate}
\end{proposition}

{\it Proof}. $2$. The proof is literally the same as the proof  to
Proposition \ref{Usl2p3}, and we leave it for  the reader. $\Box $

We left it for the reader to prove that (using a proper embedding
into  the division ring of the type $\CE_{q'}$)  the {\em quantum
Heisenberg} algebra \cite{KS}, \cite{Ma}:
$$\CH_q =K\langle X,Y,H \,|\,XH=q^2HX,\,YH=q^{-2}HY,\,XY-q^{-2}YX=q^{-1}H\rangle $$
and the  {\em Witten's first deformation}:
$$E=K\langle E_0, E_-, E_+\, | \,  [E_0,E_+]_q=E_+,
\,\,[E_-,E_0]_q=E_-,\,\, [E_+,E_-]=E_0-(q-1/q)E^2_0\rangle$$ admit
the Dixmier partition and satisfy  both the homogeneous
centralizer condition and the commutative centralizer condition
 (where $q\in K^*$ is not a root of $1$, $[a,b]_q:= qab-q^{-1}ba)$.


\section{The Ring $\CD (X)$ of Differential Operators on a Smooth
Irreducible Algebraic Curve $X$ }\label{CDXcc3}

In this section, let $K$ be an {\em algebraically closed } field
of characteristic zero, let $X$ be a {\em smooth irreducible
algebraic curve} over $K$. The coordinate algebra $\OO =\OO (X)$
 on $X$ is a {\em finitely generated regular domain of
Krull dimension} $1$. We denote by $Q=Q(X)$ its quotient field.
Let $\CD (X)=\CD (\OO (X))$ be the ring of differential operators
on $X$. Recall the definition and basic properties of the algebra
$\CD (\OO )$ (see \cite{MR}, Ch. 15, for details).

A ring of ($K$-linear) {\em differential operators} $\CD (\OO
)=\cup_{i\geq 0}\CD (\OO )_i$ on $X$ is a subalgebra of ${\rm
End}_K(\OO )$ where  $\CD (\OO)_0=\{ u\in \End_K(\OO ): \,
ur-ru=0, \; {\rm for \; all}\; r\in \OO\}=\End_{\OO}(\OO )\simeq
\OO $,
$$ \CD (\OO )_i=\{ u\in \End_K(\OO ):\, ur-ru\in \CD (\OO )_{i-1},\; {\rm for \; all \; }\;
 r\in \OO\}.$$
$\CD (\OO )_i\CD (\OO )_j\subseteq \CD (\OO )_{i+j}$ for all
$i,j\geq 0$, so $\CD (\OO )$ is a positively filtered algebra. We
say that an element $u\in \CD (\OO )_i\backslash \CD (\OO )_{i-1}$
has {\em order} $i$ denoted $\ord (u)$. The subalgebra $\D (\OO )$
of $\CD (\OO )$ generated by the coordinate  algebra $\OO \equiv
\End_{\OO}(\OO )$ and the set ${\rm Der}_K (\OO )$ of all
$K$-derivations of the algebra $\OO $ is called the {\em
derivation ring} of $\OO $.

\begin{itemize}
\item $\Der_K(\OO )$ {\em is a finitely generated projective $\OO
$-module of rank $1$.}
\item  $\CD (\OO )=\D (\OO ) $.
\item  $\CD (\OO )$ {\em is a central simple (left and right) Noetherian domain
of Gelfand-Kirillov dimension $2$.}
\item {\em If} $S$ {\em is a multiplicatively closed subset of} $\OO $ {\em then}
 $S$ {\em is a (left and right) Ore set of} $\CD (\OO ) $, {\em and}
 $\CD (\S1 \OO)=\S1 \CD (\OO ) $.
\item {\em The associative graded algebra} $\gr \, \CD (\OO ) =
\bigoplus_{i\geq 0}  \CD (\OO )_i/\CD (\OO )_{i-1}$ {\em is a
commutative domain, and so $\ord (uv)=\ord (u) +\ord (v)$ for
nonzero elements $u,v\in \CD (\OO )$.}
\item {\em If $\OO =K[x]$ is a polynomial algebra then $\CD
(K[x])=A_1$ is the Weyl algebra.}
\item {\em Let $\m $ be a maximal ideal of $\OO $. Then there
exists an element  $c=c(\m )\in \OO\backslash \m $ such that the
localization $\CD (\OO )_c$ of $\CD (\OO )$ at the powers of the
element $c$ is an Ore extension   $\OO_c[t;\d ]$  of $\OO_c$,
where $\d \in \Der_K(\OO_c)$, and $\d =\frac{d}{dx}$ where $x\in
Q(X)$ is a transcendence basis for the field $Q(X)$. Clearly,
$\ker\, \d =K$ in $Q(X)$ (since $K=\bK $).}
\end{itemize}

Let $R$ be a ring, and let $\d $ be a derivation of the ring $R$.
The {\em skew polynomial ring} (or the {\em Ore extension})
$T=R[t;\d ]$ is a ring generated freely over $R$ by an element $t$
subject to the defining relation $tr=rt+\d (r)$ for all $r\in R$.
An element $a$ of $T$ is a unique sum $\sum a_it^i$ where $a_i\in
R$.   $\deg (a) =\max \{ i\, | \, a_i\neq 0\}$ is called the {\em
degree} of $a$, $\deg (0):=-\infty $. A nonzero element $a$ of
degree $n$ can be written  as $a_nt^n+\cdots $ where by three dots
we denote the  smaller terms, $a_nt^n$ and $a_n$ are called the
{\em leading term} and the {\em leading coefficient} of $a$
respectively. If $R$  is a domain then $\deg (ab)=\deg (a) +\deg
(b)$ for all $a,b\in T$, and so $T$ is a domain.

 The algebra $\CD (X )$ is a Noetherian domain, so its (left
and right) {\em quotient ring} $\G (X)$ is a {\em division ring}.
 We fix the element $c$. Then we have the chain of algebras
\begin{equation}\label{chDcG}
 \CD (X )\subseteq \CD (X )_c=\OO_c[t;\d ]\subseteq Q(X)[t;\d
 ]\subseteq\G (X).
 \end{equation}

Our first goal is to  prove that the algebras $\CD (X)$  and $\G
(X)$ satisfy the commutative centralizer condition (Corollary
\ref{CDGccc}). We use the same strategy as before:   we embed $\G
(X)$ into a bigger division ring $\CR $ that satisfies the $ccc$
($\CR$ will be a completion of the Ore extension $Q(X)[t; \d ]$
with respect to the (additive) valuation determined by  the degree
function).

The centre of the algebra $\CD (X)$  is $K$, and the
Gelfand-Kirillov dimension is $\GK (\CD (X))=2$, so, by Corollary
\ref{NorD}, the ring of differential operators $\CD (X)$ admits
the Dixmier partition:
$$\CD (X)\backslash K=\cup_{i=1}^5\D_i(\CD (X)).$$

By the very definition of the  ring of differential operators $\CD
(X)$, $\OO (X)\backslash K\subseteq \D_1(\CD (X))$.

Let $L$ be a ring, and let $\d $ be a derivation of $L$. The {\em
formal pseudo-differential operator ring} $\CR :=L((t^{-1};\d ))$
consists of {\em inverse Laurent } series
$a=\sum_{i=-\infty}^na_it^i$ with $n\in \mathbb{Z}$ and
coefficients $a_i\in L$ where
$$ tr=rt+\d (r)\;\; {\rm for \; all}\;\; r\in L.$$
If $a\neq 0$ then the {\em degree} $\deg (a)$ of $a$  is the
maximal integer $n$ such that $a_n\neq 0$, and so $a=a_nt^n+\cdots
$ where $a_nt^n$ (resp. $a_n$) is the {\em leading term} (resp.
the {\em leading coefficient}) of $a$, and three dots denote terms
of smaller degree. We define $\deg (0)=-\infty $. Suppose that $L$
is a domain and $a,b\in \CR $. Then $\deg (ab)=\deg (a)+\deg (b)$
and $\deg (a+b)\leq \max \{ \deg (a),\deg (b)\}$. The subset
$\CR_-$ of $\CR $ that consists  of elements degree $\leq 0$ is a
subring of $\CR $, and $\CR_-t^{-1}$ is a (two-sided) ideal of
$\CR_-$ such that $\CR_-/\CR_-t^{-1}\simeq L$. If $L$ is a
division ring then $\CR $ is a division ring. The skew polynomial
ring $L[t;\d ]$ is a subring of $\CR $.

The next result is a version of the result of K. R. Goodearl
\cite{goodearlcen}, Theorem 3.5.
\begin{lemma}\label{DSccr}
Let $L$ be a proper field extension of the field $K$, let $\d $ be
a derivation of the field $L$ such that $\ker \, \d =K$. Then the
formal pseudo-differential operator ring $\CR =L((t^{-1};\d ))$
 satisfies the commutative centralizer condition.
\begin{enumerate}
\item If $u\in L\backslash K$, then $C(u)=L$.
\item If $u\in \CR \backslash L$, $\deg \, u=0$, and the leading coefficient
$u_0$ of $u$ does not belong to $K$ then $C(u)\subseteq \CR_-$,
and the map $C(u)\ra L=\CR_-/\CR_-t^{-1}$, $c\ra c+\CR_-t^{-1}$,
is a $K$-algebra isomorphism, $[C(u):K(u)]=[L:K(u_0)]<\infty $.
\item If $u\in \CR \backslash L$, and  $\deg \, u\neq 0$, then
 $C(u)=K[[v, v^{-1}]]$ is a Laurent series field  where $v$ is an element of $C(u)$ which has
 least positive degree.  $[\, C(u):K[[u, u^{-1}]] \, ]=|{\deg \, u\over \deg \, v}|<\infty $.
\end{enumerate}
\end{lemma}

{\it Proof}. Note that $\CR $ is a division ring, and so the
centralizer $C=C(u, \CR )$ of any nonzero element $u$ of $\CR $ is
a
 division ring.

$1$. Clearly, $L\subseteq C$. Suppose that $L\neq C$. Then we can
choose a nonzero element, say $a$, of $C$ of nonzero degree, say
$n$. Let $\alpha t^n$ $(\alpha  \in L)$ be its leading term. The
leading term of the commutator $[a,u]=0$ is  equal to $n\alpha
t^{n-1}\d (u)\neq 0$ since $u\not\in K=\ker \, \d $, a
contradiction. So, $C=L$.

$2$. By the assumption, $u=\sum_{j\leq 0}\, u_jt^j$ and $u_0\in
L\backslash K$. Let $0\neq a=\sum \, a_it^i \in C$. By the
previous argument, the leading term of $a$, say $a_dt^d$, belongs
to $L$, thus $\deg\, a=d=0$,   $a=\sum_{i\leq 0} \, a_it^i $, and
$C\subseteq \CR_-$. Now,
$$ 0=au-ua=\sum_{s\leq -1} (s\d (u_0)a_s+b_s)t^{s-1}$$
where the element $b_s\in L$ depends only on $a_0, \ldots ,
a_{s-1}$, and $u$. Equating the coefficients of $t^s$ to $0$ in
the identity above we obtain the infinite
 system of equations $s\d (u_0)a_s+b_s=0$, $s\leq -1$,
 with unknowns $a_s$. For a given $a_0$, the system has a unique
solution since $s\d (u_0)\neq 0$ $(u_0\not\in K=\ker\, \d$, ${\rm
char}\, K=0)$ and $L$ is a field.  So, the map $C\ra L$, $c\ra
c+\CR_-t^{-1}$, is a $K$-algebra isomorphism, hence
 $[C(u):K(u)]=[L:K(u_0)]$.

$3$. We can write $u$ as $\alpha t^n+\cdots $ for some $0\neq
\alpha \in L$ where $n:=\deg \,  u\neq 0$. Let $0\neq  b\in C(u)$.
 Then $b=\beta t^m+\cdots $ for some $0\neq \beta \in L$ where
$m:=\deg \, b$. The elements $u$ and $b$ commute, so
$$ 0=[u,b]=[\alpha t^n+\cdots , \beta t^m+\cdots ]=(n\alpha \d (\beta
)-m\beta \d (\alpha ))t^{n+m-1}+\cdots = \d
(\frac{\beta^n}{\alpha^m})\frac{\alpha^{m+1}}{\beta^{n-1}}t^{n+m-1}+\cdots
, $$
 and so $\frac{\beta^n}{\alpha^m}\in \ker\, \d =K$. The field $K$
 is an algebraically closed field, so, for an  element $\beta \in
 C(u)$  of fixed degree $m$, its leading coefficient
$\beta $ is uniquely determined by $m$ up to a scalar factor of
$K^*$.
 The set
$$ G:=\{ \deg \, c \, | \,  \; 0\neq c\in C(u)\}$$
is a nonzero subgroup of $\mathbb{Z}$ that contains $\mathbb{Z}n$,
 hence $G=\mathbb{Z}l$ where $l$ is the least positive
 element of $G$ (and so $l|n$). Fix an element $v\in C(u)$ with $\deg \, v=l$. We
 have proved above that, for a given $\deg\, c$, the leading term
 of $c\in C(u)$ is  unique (up to a factor of $K^*$). Then, by
 the choice of $v$, we can choose $\l_1\in K^*$ such that the
 elements $\l_1v^{i_1}$ and $b$  have the same leading term where
 $i_1:=nl^{-1}$. Hence  the element
 $b_1:=b-\l_1v^{i_1}$ has degree $\deg \, b_1<\deg \, b$. Repeating the same argument
 for the element $b_1\in C(u)$, we can find $\l_2\in K$ and
 $i_2 \in \mathbb{Z}$ such that the element $b_2:=b_1-\l_2v^{i_2}\in C(u)$
 has  degree $\deg \, b_2<\deg \, b_1$. Continuing in this way, we can find infinitely
many elements of $C(u)$,
$$b_\nu :=b_{\nu -1}-\l_\nu v^{i_\nu } , \; \l_\nu \in K, \; \nu \geq 1,$$
such that $\deg \, b_1>\deg \, b_2>\cdots $ and $i_1>i_2 >\cdots
$.
 So,
$$ b=\sum_{\nu \geq 1}\, \l_\nu v^{i_\nu } \in K[[v,v^{-1}]],$$
thus $C(u)\subseteq K[[v,v^{-1}]]$. The reverse inclusion is
evident since the elements $u$ and $v$ commute, thus
  $C(u)= K[[v,v^{-1}]]$. A degree argument shows that
 $[C(u):K[[u, u^{-1}]] ]=|{\deg \, u\over \deg \, v}|<\infty $.
 $\Box $

\begin{corollary}\label{CDGccc}
 Let $\CD (X)$ be the ring of differential operators on a smooth irreducible
algebraic curve $X$, and let $\G (X)$ be its quotient division
ring. Then $\CD (X)$ and $\G (X)$ satisfy the commutative
centralizer condition.
\end{corollary}

{\it Proof}. By (\ref{chDcG}), the division ring $\G (X)$ is the
quotient ring for the algebra $Q(X)[t;\d ]$ which is  a subring of
the formal pseudo-differential ring $Q(X)((t^{-1};\d ))$. The
latter  is a division ring. By  the universality of localization,
the division ring   of $Q(X)[t;\d ]$, that is $\G (X)$, is a
subring of $Q(X)((t^{-1};\d ))$. So,  $\CD (X)\subseteq \G
(X)\subseteq Q(X)((t^{-1};\d ))$ and $\ker \, \d =K$ in $Q(X)$. By
Lemma \ref{DSccr}, the formal pseudo-differential operator ring
$Q(X)((t^{-1};\d ))$ satisfies the commutative centralizer
condition, hence $\CD (X)$ and $\G (X)$ satisfy the commutative
centralizer condition.  $\Box $

\begin{corollary}\label{AsubDX}
The ring  of differential operators $\CD (X)$ on a smooth
irreducible algebraic curve $X$, or any
  non-commutative subalgebra  $A$ of the division ring $\G (X)$
 with $\GK (A)<3$ admits the Dixmier partition and
 satisfies both the homogeneous centralizer condition and the
 commutative  centralizer condition.
\end{corollary}

{\it Proof}. Since $K=\bK $, $A$ is a noncommutative domain, and
$\GK (A)<3$, we see that the hypothesis of Corollary \ref{NorD}
holds, so $A$ admits the Dixmier partition. The algebra $\G (X)$
satisfies the commutative  centralizer condition, so, by Theorem
\ref{G3Apr}, the algebra $A$ satisfies both the homogeneous
centralizer condition and the commutative  centralizer condition.
$\Box $

\begin{corollary}\label{CDXDP}
{\sc (An Analog of Dixmier's Problem 6 for $\CD (X)$)} Let $\CD
(X)$ be the ring of differential operators on a smooth irreducible
algebraic curve $X$, and let $f(t)\in K[t]$ be an arbitrary
polynomial  of degree $>1$.
\begin{enumerate}
\item Let  $i=1,2,5$. If $a\in \D_i(\CD (X))$ then $f(a)\in \D_i(\CD (X))$.
\item Let $i=3,4$. If $a\in \D_i(\CD (X))$  then $f(a)\in \D_5(\CD (X))$.
\end{enumerate}
\end{corollary}

{\it Proof}. Note that if  $a\in \CD (X)\backslash K$ then
$f(a)\in \CD (X)\backslash K$. For the element $a\in \OO
(X)\backslash K$ this fact is evident since $K=\bK $, and for
$a\not\in \OO (X)$ this follows from $\ord (f(a))=\deg_t(f)\cdot
\ord (a)>0$. Clearly, the non-scalar elements  $a$ and $f(a)$
commute, so $C(a)=C(f(a))$ (Corollary \ref{PACf=Cg}.(1)) since
$\CD (X)$ satisfies the commutative centralizer condition
(Corollary \ref{CDGccc}).

For $i=1,2$, the result follows from  Proposition \ref{NAappb}.

Suppose that $a\in \D_{3,4,5}(\CD (X))$. By Theorem
\ref{CDXnonex}, $f(a)\not\in \D_{3,4}(\CD (X))$, and, by
Proposition \ref{NAappb}, $f(a)\not\in \D_{1,2}(\CD (X))$, hence
$f(a)\in \D_5(\CD (X))$. $\Box $

{\bf Question 1}. {\em Find generators for the group $\Aut_K(\CD
(X))$ of algebra automorphisms.}

{\bf Question 2}. {\em Classify (up to the action of the group
$\Aut_K(\CD (X)))$ nilpotent and semi-simple elements of $\CD
(X)$.}

{\bf Question 3}. {\em For a semi-simple element $h$  of $\CD
(X)$, is the set $\Ev (h, \CD (X))$ of eigenvalues for the inner
derivation $\ad \, h$ of $\CD (X)$  equal to $ \mathbb{Z}\rho $
for some $\rho \in K$?}

{\bf Question 4}. {\em Let $h$ and $a$  be a semi-simple and
nilpotent element of the algebra $\CD (X)$ respectively. Is the
eigenvalue algebra $D(h, \CD (X))$ and the nil-algebra $N(a,\CD
(X))$ finitely generated (Noetherian)?}

{\bf Question 5}. {\em For which $X$, every algebra endomorphism
of $\CD (X)$ is an algebra automorphism?}

\bigskip
\begin{center}
{\sc Acknowledgment }
\end{center}
\bigskip
\noindent The author would like to thank J. Dixmier for comments
\cite{Dixcom}
 on his problems from \cite{Dix}, and D. Jordan and T. Lenagan for
 pointing out on the paper of K. R. Goodearl \cite{goodearlcen}.

Department of Pure Mathematics

University of  Sheffield

Hicks Building

Sheffield S3~7RH

email: v.bavula@sheffield.ac.uk


\begin{thebibliography}{99}

\bibitem{Ami} S. A. Amitsur,  Commutative linear differential operators,
 {\it Pacific J. Math.} {\bf 8} (1958), 1--10.

\bibitem{Art-Cohnqpl} V. A. Artamonov and P. M.  Cohn,  The skew field of rational
functions on the
   quantum plane. {\em Algebra, 11. J. Math. Sci. (New York)} {\bf  93} (1999),
    no. 6, 824--829.

\bibitem{BCW}    H. Bass, E. H. Connel and D. Wright,
 The Jacobian Conjecture: reduction of degree and formal expansion of the inverse,
{\it Bull. Amer. Math. Soc. (New Series)}, {\bf 7} (1982),
287--330.

\bibitem{BassM}    H. Bass, G. Meisters, Polynomial flows in the plane,
 Adv. in Math.  {\bf 55} (1985), no. 2, 173--208.

\bibitem{bavaia}  V. V. Bavula, Generalized Weyl algebras and their representations.
(Russian) {\it Algebra i Analiz} {\bf 4} (1992), no. 1, 75--97;
translation in {\it St. Petersburg Math. J.} {\bf 4} (1993), no.
1, 71--92.


\bibitem{Bavrenqu} V. V. Bavula, A Question of Rentschler and the
Problem of Dixmier, {\it Ann. of Math.} {\bf 154} (2001), no. 3,
683--702.

\bibitem{BavDP5} V. V. Bavula, Dixmier's Problem 5 for the Weyl algebra, {\it }
(submitted).

\bibitem{BGR} W.  Borho, P.  Gabriel, and R. Rentschler,
Primideale in Einhüllenden auflösbarer Lie-Algebren (Beschreibung
durch Bahnenräume). (German) Lecture Notes in Mathematics, Vol.
357. Springer-Verlag, Berlin-New York, 1973.

\bibitem{CBH}  W. Crawley-Boevey and M.  Holland.  Noncommutative deformations of
Kleinian singularities. {\it Duke Math. J.} {\bf 92} (1998), no.
3, 605--635.

\bibitem{Dix3} J. Dixmier, Repr$\acute{e}$sentation irr$\acute{e}$ducibles des alg\'ebres de Lie resolubles,
{\it J. Math. pures et appl.} {\bf 45} (1966), 1--66.


\bibitem{Dix}
J. Dixmier,  Sur les alg\`{e}bres de Weyl, {\it Bull. Soc. Math.
France} {\bf 96} (1968), 209--242.




\bibitem{Dixcom}
J. Dixmier, a letter.


\bibitem{hodges} T. Hodges, Noncommutative deformations of type-$A$ Kleinian singularities.
{\it J. Algebra} {\bf 161} (1993), no. 2, 271--290.


\bibitem{goodearlcen} K. R. Goodearl, Centralizers in differential,
pseudodifferential, and fractional  differential operator rings.
{\it Rocky Mountain J. Math.} {\bf  13}
   (1983), no. 4, 573--618.


\bibitem{josclA1} A. Joseph, The Weyl algebra---semisimple and nilpotent
elements, {\it Amer. J. Math.} {\bf 97} (1975), no. 3, 597--615.

\bibitem{KL} G. Krause and T.  Lenagan, {\em Growth of algebras and Gelfand-Kirillov
 dimension}. Revised edition. Graduate Studies in Mathematics, 22.
American Mathematical Society, Providence, RI, 2000.


\bibitem{KS} E. E. Kirkman and L. W. Small, $q$-analogs of harmonic oscillators
and related rings, Preprint, Wake Forest University and University of California, 1992.

\bibitem{Ma} M.-P. Malliavin, L'alg\`{e}bre d'Heisenberg quantique,
 C. R. Acad. Sci. Paris, S$\acute{e}$r. 1, 317 (1993), 1099--1102.

\bibitem{Mazccc} V. Mazorchuk, A note on centralizers in $q$-deformed Heisenberg
   algebras. {\it AMA Algebra Montp. Announc.}  (2001), Paper 2, 6 pp.

\bibitem{MRext} J. C. McConnell and J. C. Robson, Homomorphisms and extensions
of modules over certain differential polynomial rings, {\it J.
Algebra} {\bf 26} (1973), 319--342.

\bibitem{MR} J. C. McConnell and J. C. Robson, {\em Noncommutative Noetherian
rings}. With the cooperation of L. W. Small. Revised edition.
Graduate Studies in Mathematics, 30. American Mathematical
Society, Providence, RI, 2001.


\end{thebibliography}
\end{document}